\documentclass[12pt,twoside,leqno]{article}
\usepackage{amsmath}
\usepackage{amssymb}
\usepackage{amsxtra}
\usepackage{amscd}
\usepackage{amsthm}
\usepackage[mathscr]{eucal}
\usepackage{xr}
\usepackage{color}

\swapnumbers

\theoremstyle{plain}

\newtheorem{sbthm}[subsubsection]{Theorem}
\newtheorem{sbprop}[subsubsection]{Proposition}
\newtheorem{sbcor}[subsubsection]{Corollary}
\newtheorem{sblem}[subsubsection]{Lemma}

\theoremstyle{definition}

\newtheorem{sbrem}[subsubsection]{Remark}
\newtheorem{sbpara}[subsubsection]{}

\newenvironment{pf}{\proof[\proofname]}{\endproof}
\newenvironment{pf*}[1]{\proof[#1]}{\endproof}

\newcommand{\N}{{\mathbb{N}}}
\newcommand{\Q}{{\mathbb{Q}}}
\newcommand{\Z}{{\mathbb{Z}}}
\newcommand{\R}{{\mathbb{R}}}
\newcommand{\C}{{\mathbb{C}}}

\newcommand{\Spec}{\operatorname{Spec}}

\newcommand{\fs}{{\rm fs \ }}
\newcommand{\gp}{\mathrm{gp}}
\newcommand{\Gm}{\mathbb{G}_m} %
\newcommand{\Gmlog}{\mathbb{G}_{m,\log}}

\newcommand{\Hom}{\operatorname{Hom}}
\newcommand{\End}{\operatorname{End}}
\newcommand{\Ext}{\operatorname{Ext}}

\newcommand{\Aut}{\operatorname{Aut}}

\newcommand{\Tr}{\operatorname{Tr}}

\newcommand{\Ker}{\operatorname{Ker}}

\renewcommand{\bar}[1]{\overline{#1}}

\def\spcheck{^\vee}

\def \et {\mathrm {\acute{e}t}}

\newcommand{\cl}{Claim}
\theoremstyle{definition}

\newtheorem*{clm*}{\cl}

\theoremstyle{plain}

\def \fs {\mathrm {fs}}
\def \ket {\mathrm {k\acute{e}t}}

\def \overc#1{\overset {\lower 0.3ex \hbox{${\;}_{\circ}$}}{#1}}

\let\refsave=\ref
\def\ref#1{\textup{\refsave{#1}}}

\newcommand{\upc}{\overset{\circ}\to}
\newcommand\Cal{\mathcal}
\newcommand\define{\newcommand}
\renewcommand\bold{\Bbb}

\define\bZ{\bold Z}
\define\bC{\bold C}
\define\bR{\bold R}
\define\bQ{\bold Q}
\define\bN{\bold N}
\define{\cS}{\Cal S}
\define{\Lie}{\mathrm{Lie}\,} %
\define{\coLie}{\mathrm{coLie}\,} %
\define{\cH}{\Cal H}
\define{\cExt}{{\Cal E}xt}
\define{\cHom}{{\Cal H}om}
\define{\cO}{\Cal O}
\define{\an}{\mathrm{an}}

\def \upcf {\overset {\lower 0.3ex \hbox{${\;}_{\circ}$}} f}
\def \upcp {\overset {\lower 0.3ex \hbox{${\;}_{\circ}$}} p}
\def \upc#1{\overset {\lower 0.3ex \hbox{${\;}_{\circ}$}}{#1}}

\newcommand{\bs}{\backslash}
\newcommand{\Sig}{\Sigma}
\newcommand{\sig}{\sigma}

\newcommand{\cA}{\mathcal A}
\newcommand{\cB}{\mathcal B}
\newcommand{\cC}{\mathcal{C}}
\newcommand{\cE}{\mathcal{E}}

\newcommand{\cG}{\mathcal{G}}

\newcommand{\cT}{\mathcal T}

\newcommand{\gr}{{\mathrm{gr}}}

\usepackage{fancyhdr}
\usepackage{amssymb}
\usepackage{graphicx}

\begin{document}
\title
{Moduli of logarithmic abelian varieties with PEL structure}
\author{Takeshi Kajiwara, Kazuya Kato, and Chikara Nakayama}
\date{}
\maketitle
\centerline{\it{Dedicated to Professor Luc Illusie}} 
\setlength{\baselineskip}{1.0\baselineskip}
\begin{abstract}
\noindent 
  We construct the fine moduli space of log abelian varieties with PEL structure, which gives a toroidal compactification of the moduli space of abelian varieties with PEL structure.
\end{abstract}

\begin{center}
{\small {\bf R\'esum\'e}}
\end{center}

\begin{quotation}
\vskip -1.5ex
{\small \noindent 
  Nous construisons l'espace des modules fin des vari\'et\'es ab\'eliennes logarithmiques de type PEL, qui donne une compactification toro\"idale de l'espace des modules des vari\'et\'es ab\'eliennes de type PEL.}
\end{quotation}

\section*{Contents}

\noindent \S\ref{s:mod}. Moduli functors

\noindent \S\ref{s:main}. The main theorems and proofs

\noindent \S\ref{s:artin}. Second proof

\noindent \S\ref{s:toroidal}. Complements 
\renewcommand{\thefootnote}{\fnsymbol{footnote}}
\footnote[0]{MSC2020: Primary 14A21; Secondary 14K10, 14D06} 
\footnote[0]{Keywords: degeneration of abelian varieties, log geometry}

\section*{Introduction}
  We construct toroidal compactifications
of the moduli space of abelian varieties with PEL structure as
the fine moduli space of log abelian varieties with PEL structure.
  Our toroidal compactification is essentially the same as the one constructed by Lan (\cite{Lan}), but a good thing of our approach is that we have a  good moduli functor. By virtue of the remarkable feature of log abelian varieties that they are some kind of compactified degenerate abelian varieties but still have group structures, our formulation of the moduli problem  is exactly parallel to the nondegenerate case. In particular, in our formulation, the endomorphism ring (the E in PEL) is really the endomorphism ring for the group structure of a log abelian variety, and this makes our theory transparent.  
(The polarization P and the level structure L are also formulated by using the group structure.)

  More precisely, we prove the following (1) and (2). 

(1) The moduli functor of $g$-dimensional log abelian varieties with PEL structure is represented by a proper and log smooth log algebraic space in the second sense (Theorem \ref{t:Main}).

(2) The moduli functor of $g$-dimensional log abelian varieties with PEL structure and with local monodromies in a prescribed admissible cone decomposition is represented by a proper and log smooth algebraic space with fs log structure 
(Theorem \ref{t:Main_Sig}).

  We started our study of log abelian varieties about thirty years ago to understand the paper \cite{Fujiwara} of K.\ Fujiwara by log geometry. This paper \cite{Fujiwara} also constructs toroidal compactifications of moduli spaces of abelian varieties with PEL structure.  
  The prototypes of all our results and techniques of proofs can be found in his paper. 

  The organization of this paper is as follows. 
  In Section \ref{s:mod}, we introduce moduli functors. 
  Section \ref{s:main} consists of main results and their proofs. 
  With the feature of log abelian varieties explained in the above, our proofs are also parallel to the nondegenerate case (cf.\ \cite{GN}), that is, we reduce to the non-coefficient case (\cite{KKN7}) by the representability of $\cHom$-sheaves. 
  In Section \ref{s:artin}, we improve log Artin's criterion (Theorem \ref{logAr}) and give alternative proofs for the main results. 
  In Section \ref{s:toroidal}, we give some complements: 
We compare our space in the above (2) with the analytic moduli space constructed in \cite{KKN1} (Proposition \ref{overC}) and also with Lan's compactification (Proposition \ref{p:comparison}).
  \smallskip

\noindent {\sc Acknowledgments.}
The authors thank B.\ C.\ Ng\^ o for advice on Shimura varieties. 
The first author is partially supported by JSPS, Kakenhi (C) No.\ 24540035, Kakenhi (C) No.\ 15K04811, and Kakenhi (C) No.\ 20K03555. 
The second author is partially  supported by NFS grants DMS 1303421, DMS 1601861, and DMS 2001182.
The third author is partially supported by JSPS, 
Kakenhi (B) No.\ 23340008, Kakenhi (C) No.\ 16K05093, and 
Kakenhi (C) No.\ 21K03199.

\section{Moduli functors}
\label{s:mod}

\subsection{$B$, $V$ and $\psi$}\label{s:BVpsi}

\begin{sbpara}\label{BV}
We fix a finite dimensional semisimple algebra $B$ over $\Q$ with a map $B\to B\;;\; x\mapsto x^*$ such that $(x+y)^*=x^*+y^*$, $(xy)^*=y^*x^*$, $(x^*)^*=x$ for all $x,y\in B$ and such that $\Tr_{B/\Q}(xx^*)>0$ for all $x\in B \smallsetminus \{0\}$. 

We fix a finitely generated $B$-module $V$ and a nondegenerate skew-symmetric $\Q$-bilinear form $\psi: V \times V \to \Q(1):=\bQ\cdot 2 \pi i$ such that $\psi(bx, y)= \psi(x, b^*y)$ for all $x,y\in V$ and for all $b\in B$. 

Let $G$ be the reductive algebraic group over $\Q$ defined as follows. For a commutative ring $R$ over $\Q$, $$G(R)= \{(a,c)\in \text{Aut}_{R\otimes B}(R\otimes V) \times R^\times\;|\; \psi(ax,ay)= c\psi(x,y) \;\text{ for all}\; x,y\}.$$

Let $2g=\dim_\Q(V)$. 

\end{sbpara}

\begin{sbpara}\label{Sdata} We assume that we are given an $\R$-algebra homomorphism $h: \bC \to \End_{B \otimes\R}(V \otimes \R)$ such that $\psi(h(z)x,y)=\psi(x,h(\bar z)y)$ for all $x,y \in V \otimes\bR$ and for all $z \in \bC$, satisfying the following condition: 

  The $\bR$-bilinear pairing 
$\frac{1}{2 \pi i}\psi(\cdot,h(i)\cdot)\colon 
(V \otimes\bR) \times (V \otimes\bR) \to \bR$
is symmetric and positive definite. 

For such an $h$, consider the homomorphism of algebraic groups over $\R$
$$h_0\;:\; S_{\C/\R}= \text{Res}_{\C/\R}({\mathbb G}_{m,\C}) \to G_\R$$
which induces $$S_{\C/\R}(\R)=\C^\times \to G_\R(\R)\subset \text{Aut}_\R(V \otimes \R) \times \R^\times\;;\; z \mapsto (h(z), (z\bar z)^{-1}).$$
Then the $G(\R)$-conjugacy class of $h_0$ 
gives Shimura data (\cite{De} 1.5). 
  Let $F_0$ be the reflex field  $E(G, h_0)$ (\cite{De} 3.7) of the Shimura data. 

\end{sbpara}

\subsection{Moduli functors without fans}\label{s:functor} 
We try to follow the notation of Lan's book (\cite{Lan}) as much as possible. 

\begin{sbpara}\label{BU2} Let $\square$ be a set of prime numbers.  Let $\Z_{(\square)}$ and $O_{F_0,(\square)}$ be the localization of $\Z$ and $O_{F_0}$, respectively, by inverting all prime numbers outside $\square$. 

Assume that we are given a  $\Z_{(\square)}$-subalgebra $\cO_{(\square)}$ of $B$ which is finitely generated as a $\Z_{(\square)}$-module such that $\Q\otimes \cO_{(\square)}= B$, 
such that $\cO_{(\square)}^*=\cO_{(\square)}$, and 
such that the center of $\cO_{(\square)}$ is finite \'etale over $\Z_{(\square)}$. 
  We fix such an $\cO_{(\square)}$. 
  Assume also that we are given a finitely generated $\cO_{(\square)}$-submodule $L_{(\square)}$ of $V$ such that $V=\Q\otimes L_{(\square)}$, such that $\psi(L_{(\square)}, L_{(\square)})\subset \Z_{(\square)}(1)$, and such that $\psi: L_{(\square)} \times L_{(\square)}\to \Z_{(\square)}(1)$ is a perfect pairing. 
  If $B$ involves a simple factor of type D (cf.\ \cite{Lan} Definition 1.2.1.15), we assume further that $\square$ is prime to $2$. 

  Then $G$ extends to a smooth affine group scheme $G_{L_{(\square)}}$ over $\bZ_{(\square)}$ defined as follows. For a commutative ring $R$ over $\bZ_{(\square)}$, 
$$G_{L_{(\square)}}(R)= \{(a,c)\in \text{Aut}_{R\otimes\cO_{(\square)}}
(R\otimes L_{(\square)}) \times R^\times\;|\; \psi(ax,ay)= c\psi(x,y) \;\text{ for all}\; x,y\}.$$
\end{sbpara}

\noindent {\it Remark.}
   Fujiwara (\cite{Fujiwara}) assumes that any $\bR$-simple factor of $G_{\bR}^{\mathrm {ad}}$ has type A or C.

\begin{sbpara}\label{FK12} 

Let $\hat \Z^{\square}=\prod_{\ell\notin \square} \Z_{\ell}$, where $\ell$ ranges over prime numbers which do not belong to $\square$, and let $\hat \Q^{\square}= \Q\otimes \hat \Z^{\square}$. 
Fix a compact open subgroup $\cH$ of $G(\hat \Q^{\square})$. 
  We assume that $\cH$ is {\it neat} in the sense that 
any element $g=(g_\ell)_{\ell \not\in\square}$ of $\cH$ satisfies $\bigcap_{\ell \not\in \square} (\bar \bQ^{\times} \cap \Gamma_{g_\ell})_{\mathrm{tors}}=\{1\}$, where $\Gamma_{g_\ell}$ is the subgroup of $\bar \bQ_\ell^{\times}$ generated by eigenvalues of $g_\ell$ (cf.\ \cite{Lan} Definition 1.4.1.8). 
\end{sbpara}

\begin{sbpara}\label{FHbarH}

We define moduli functors $$F_\cH : (\text{sch}/O_{F_0,(\square)})\to (\mathrm{set})$$
and 
$${\bar F}_\cH : (\text{fs}/O_{F_0,(\square)})\to (\mathrm{set})$$
as in \ref{FK22} below. 
  Here $(\text{sch}/O_{F_0,(\square)})$ (resp.\ $(\text{fs}/O_{F_0,(\square)})$) means the category of schemes (resp.\ fs log schemes) over $O_{F_0,(\square)}$. 

  The former is represented by a smooth scheme 
over $O_{F_0,(\square)}$ (cf.\ \cite{Lan} Theorem 1.4.1.11).
  We will show that the latter is a log smooth log algebraic space in the second sense (Theorem \ref{t:Main}).
\end{sbpara}

\begin{sbpara}\label{FK22} 
  Let $F_\cH(S)\; (\text{resp.}\; \bar F_\cH(S))= \{(A, \iota, p, \eta)\}$.
This $4$-ple is considered in the category of abelian schemes (resp.\ log abelian varieties) in which Hom is $\Hom(A, A') \otimes \Z_{(\square)}$ in the usual sense. 

0. $A$ is an abelian scheme (resp.\ a log abelian variety) over $S$ of dimension $g$, considered modulo isogeny of degree prime to $\square$ (that is, in the category where $\Hom$ is defined as  $\Hom\otimes \Z_{(\square)}$). 

1. $\iota$ is a homomorphism $\cO_{(\square)}\to \text{End}(A) \otimes \Z_{(\square)}$ satisfying a certain condition described in \ref{cond1} below. 

2. $p$ is an element of $\Hom(A, \cE xt(A, \Gmlog)) \otimes \Z_{(\square)}$ written as $n^{-1}p'$ for some integer $n$ which is prime to $\square$ and for some  polarization $p': A\to \cE xt(A, \Gmlog)$ whose degree is prime to $\square$ satisfying a certain condition described in \ref{cond2} below. 
  (Recall that a polarization on a log abelian variety $A$ is regarded as a homomorphism $A \to \cExt(A,\Gmlog)$ satisfying a certain condition (\cite{KKN5}).)

3. $\eta$ is a class of a $B\otimes \hat \Q^{\square}$-isomorphism modulo $\cH$ on the pro-\'etale site (resp.\ pro-k\'et site) of $S$ from $ V^{\square}A:= (\prod_{\ell \notin \square} T_{\ell}A)\otimes \Q$ to  $V \otimes  \hat \Q^{\square}$ satisfying a certain condition described in \ref{cond3} below, where $T_{\ell}A$ is the Tate module of $A$. 
  (Here it is important that $\Ker(\ell^m: A \to A)$ $(m>0)$ is a locally constant sheaf on the kummer log \'etale (= k\'et) site by \cite{KKN4} Proposition 18.1 (3).)
\end{sbpara}

\noindent 
{\it Remark.}  
  Note that in the case of $F_\cH$, the above definition of level structure in the above 3 is equivalent to the one in \cite{Lan} 1.3.
  We can write the above 3 for $\bar F_\cH$ with the usual k\'et site in the same way as \cite{Lan} uses the usual \'etale site. 
  We prefer the pro-k\'et site to simplify the description. 
  Here the pro-k\'et site is defined as the site of pro-objects of k\'et objects as in the non-log case. 

\begin{sbpara}\label{cond1}
  The condition in \ref{FK22} 1 is as follows. 
  Let $V_0$ be the subspace of $V \otimes \bC$ on which $h(z)$ acts by the multiplication by $z$ ($z \in \bC$). 
  The condition is that $\det_{\cO_{(\square)}\vert\Lie_{A/S}}$ agrees with the image of $\det_{\cO_{(\square)} \vert V_0}$ under the structural homomorphism $O_{F_0,(\square)} \to \cO_S(S)$.
  Here, $\det_{\cO_{(\square)}\vert\Lie_{A/S}}$ and 
$\det_{\cO_{(\square)}\vert V_0}$ are defined as in \cite{Lan} Definition 1.1.2.21 as follows. 
  Take a $\bZ_{(\square)}$-basis $\alpha_1,\ldots, \alpha_t$ of $\cO_{(\square)}$. 
  Let $\alpha_1\spcheck,\ldots,\alpha_t\spcheck$ be the dual basis of $\cO_{(\square)}\spcheck:=\Hom_{\bZ_{(\square)}}(\cO_{(\square)},\bZ_{(\square)})$. 
  Consider the action of $\cO_{(\square)}$ on $\Lie_{A/S}$ via $\iota$.  
  Let $\det_{\cO_{(\square)}\vert\Lie_{A/S}}$ be the image of 
the determinant of $X_1 \alpha_1+\cdots +X_t\alpha_t$ on $\Lie_{A/S}$ 
by the homomorphism 
$\cO_S(S)[X_1,\ldots,X_t] \to \cO_S(S)[\cO_{(\square)}\spcheck]$ sending $X_j$ to $\alpha_j\spcheck$ ($1 \leq j \leq t)$.
  Next, $\cO_{(\square)}$ acts on $V_0$ via the action of $B$ on $V \otimes \bC$.
  By using it, $\det_{\cO_{(\square)}\vert V_0}$ is defined by the image of 
the determinant of $X_1\alpha_1+\cdots +X_t\alpha_t$ on $V_0$ by the homomorphism 
$\bC[X_1,\ldots,X_t] \to
\bC[\cO_{(\square)}\spcheck]$ sending $X_j$ to $\alpha_j\spcheck$ ($1 \leq j \leq t)$, and belongs to $O_{F_0,(\square)}[\cO_{(\square)}\spcheck]$. 
 \end{sbpara}

 \begin{sbpara}\label{cond2}

The condition in \ref{FK22} 2 is that for $b\in \cO_{(\square)}$, the following diagram is commutative in the category modulo isogeny of degree prime to $\square$: 
$$\begin{matrix}  A  & \overset{b^*}\to & A \\
{}^p \downarrow {\;}&& {\;}\downarrow {}^p\\
\cE xt(A, \Gmlog) & \overset{\text{by $b$}}\to & \cExt(A, \Gmlog).\end{matrix}$$
\end{sbpara}

\begin{sbpara}\label{cond3}

The condition in \ref{FK22} 3 is that $V^{\square}A \times V^{\square}A\overset{1\times p}\to V^{\square}A \times V^{\square}\cExt(A,\Gmlog)\to \hat \Q^{\square}(1)$ is compatible with $\psi$, 
where the last arrow is the Weil pairing (cf.\ \cite{KKN7} 9.1). 

\end{sbpara}

\subsection{Moduli functors with fans}\label{s:fan}

Let the notation and the assumptions be as in Section \ref{s:functor}. We consider a compatible family $\Sig$ of complete fans and we define a subfunctor $\bar F_{\cH,\Sig}$ of $\bar F_\cH$. 

\begin{sbpara}
For a positive semidefinite symmetric bilinear form $\mu$ on an $\R$-vector space $Z$, we call the subspace $\Ker(\mu):=\{z\in Z\;|\; \mu(y, z)=0\text{ for all }y\in Z\}$ the {\it kernel} of $\mu$. The kernel of $\mu$ coincides with $\{z\in Z\;|\; \mu(z, z)=0\}$.  
\end{sbpara}

\begin{sbpara}\label{cofan} Let $Z$ be a finitely generated $B$-module.  Let $C^+(Z)$ be the set of all positive semidefinite symmetric $\R$-bilinear forms $\mu$ on $Z \otimes\R$ with rational kernel such that $\mu(bx, y)= \mu(x, b^*y)$ for all $x,y\in Z\otimes \R$ and $b\in B$. 

By a {\it complete fan} for $Z$, we mean a fan of finitely generated rational cones  in   $C^+(Z)$ satisfying the following conditions (i) and (ii).

(i)  The support of this fan is $C^+(Z)$.  

(ii) There are a finite subset $S$ of this fan and an arithmetic subgroup $\Gamma$ of $\Aut_B(Z)$ such that all elements of this fan are written as $\gamma \sig$ for $\gamma\in \Gamma$ and $\sig\in S$. (Equivalently, for each arithmetic subgroup $\Gamma$ of $\Aut_B(Z)$,  there is a finite subset $S$ of this fan such that  all elements of this fan are written as $\gamma \sig$ for $\gamma\in \Gamma$ and $\sig\in S$.)

\end{sbpara}

\begin{sbpara}\label{Sig} A {\it compatible family of complete fans} for $(B, V, \psi,\cH)$ is to give a complete fan $\Sig(W,g)$ for $W\otimes \Q$ (\ref{cofan})  for each pair $(W, g)$ of a finitely generated $\cO_{(\square)}$-module  $W$ which is torsion-free as a $\Z$-module and a surjective $B\otimes \hat \Q^{\square}$-homomorphism $$g: V \otimes \hat \Q^{\square}\to W\otimes \hat \Q^{\square}$$ whose kernel $J$ satisfies
$J^{\perp}\subset J$ (that is, $\psi(J^{\perp}, J^{\perp})=0$\;;\; here $J^{\perp}\subset V\otimes\hat \Q^{\square}$ denotes the annihilator of $J$ under $\psi$), satisfying the following conditions (i) and (ii).

(i) $\Sig(W, gk)= \Sig(W, g)$ for $k\in \cH$.

(ii) Let $W'$ be a finitely generated $\cO_{(\square)}$-module which is torsion-free as a $\Z$-module and  let $\gamma: W\to W'$ be a surjective $\cO_{(\square)}$-homomorphism. 
  Then $\Sig(W', \gamma g)= \gamma \Sig(W, g)$. Here $\gamma g$ denotes the composite homomorphism $$V\otimes \hat \Q^{\square} \overset{g}\to W\otimes \hat \Q^{\square}\overset{\gamma}\to W'\otimes \hat \Q^{\square},$$ 
and   $\gamma \Sig(W,g)$ is the set of cones $\sig$ in $C^+(W'\otimes \Q)$ such that the image of $\sig$ in $C^+(W\otimes \Q)$ under the map $C^+(W'\otimes \Q) \to C^+(W\otimes \Q)$ induced by $\gamma$ belongs to $\Sig(W,g)$. 
\end{sbpara}

\begin{sbprop}[\cite{AMRT} Chapter II, \cite{Lan} Proposition 6.3.3.5]
A compatible family of complete fans for $(B, V, \psi,\cH)$ exists. 
\end{sbprop}

\begin{sbpara}\label{FKSig2}  Let $\Sig=(\Sig(W,g))_{(W,g)}$ be a compatible family of complete fans for $(B, V,\psi, \cH)$.

We define a moduli functor $$\bar F_\cH\supset \bar F_{\cH, \Sig}: (\fs/O_{F_0, (\square)})\to (\mathrm{set})$$
as below.  We will show that $\bar F_{\cH, \Sig}$  is represented by an algebraic space with an fs log structure (Theorem \ref{t:Main_Sig}).

\end{sbpara}

\begin{sbpara}\label{cond4-0} The following (1)--(3) are  preparations   for the definition of the moduli functor $\bar F_{\cH, \Sig}$ with fan.

(1) For a log abelian variety $A$ over an fs log scheme $S$, we have a  sheaf of abelian groups $\bar Y$ on the big \'etale site $(\fs/S)_{\et}$ of $S$ whose all stalks are finitely generated free abelian groups (\cite{KKN2} 4.1, 9.3). 
A polarization of  $A$ induces a $\Z$-bilinear form $\bar Y \times \bar Y\to (\Gmlog/\Gm)_S$ which sends $(y, y)$ for a local section $y$ of $\bar Y$ to $(M/\Gm)_S$ and which is nondegenerate in the sense that if $(y, y)$ is sent to $1\in (M/\Gm)_S$ then $y=0$ (\cite{KKN2} 9.1). 
  Here $M$ is the sheaf $U \mapsto \Gamma(U,M_U)$. 
If $\ell$ is a prime number which is invertible on $S$, we have a canonical surjective homomorphism $T_{\ell}A \to \bar Y \otimes \Z_{\ell}$ (\cite{KKN4}  18.9.1). This 
$\bar Y \otimes \Z_{\ell}$ is the (weight $0$)-quotient of $T_{\ell}A$. For a polarization $p$ and for the pairing $T_{\ell}A \times T_{\ell}A \to \Z_{\ell}(1)$ induced by $p$, if we denote by $J$ a stalk of the kernel of  $T_{\ell}A \to Y \otimes  \Z_{\ell}$ (this $J$ is the (weight $\leq -1$)-part of the stalk of $T_{\ell}A$), the annihilator of $J$ in $T_{\ell}A$ under this pairing (it is the (weight $-2$)-part of the stalk of $T_{\ell}A$) is contained  in $J$.  

(2) Let $A$ and $S$ be as in (1) and assume that the underlying scheme of $S$ is the spectrum of a separably closed field. Denote $S$ by $s$ here.   
The log fundamental group $\pi_1^{\log}(s)=\Hom((M^{\gp}_s/\cO^\times_s)_s, \prod_{\ell} \Z_{\ell}(1))$, where $\ell$ ranges over all prime numbers which are invertible on $s$, acts on $V^{\square}A$. The action is described as follows (\cite{KKN4} Proposition 18.11). For $\gamma\in \pi_1^{\log}(s)$, the action of $\gamma$ on $V^{\square}A$ satisfies $(\gamma-1)^2=0$, and 
$\log(\gamma)=\gamma-1$ coincides with the composition
$$V^{\square}A \to \bar Y_s \otimes \hat \Q^{\square} \overset{\gamma}\to \Hom(\bar Y_s, \hat \Q^{\square}(1)) \to \Hom_{\hat \Q^{\square}}(V^{\square}A, \hat \Q^{\square}(1))\cong V^{\square}A$$ 
in which the arrow $\gamma$ is induced from $\bar Y_s \otimes \bar Y_s \to (M^{\gp}_s/\cO^\times_s)_s\otimes \Q$ and $\gamma$. 
  This action of $\pi_1^{\log}(s)$ respects the pairing $V^{\square}A \times V^{\square}A \to \hat \Q^{\square}(1)$ induced by the polarization. 

(3) For $(W, g)$ as in \ref{Sig}, a symmetric bilinear form $\mu: W \times W \to \hat \Q^{\square}$ gives a $\hat \Q^{\square}$-linear map
$$N_\mu: V \otimes \hat\Q^{\square}\to V\otimes \hat \Q^{\square}(-1)$$ characterized by 
$$\psi(N_\mu x, y) = \mu(\bar x, \bar y)$$
for all $x,y\in V\otimes \hat \Q^{\square}$, where $\bar x, \bar y$ denote the canonical images of $x, y$ in $W \otimes \hat \Q^{\square}$. This $N_\mu$ satisfies $N_\mu^2=0$ and
$$\psi(N_\mu x, y)+\psi(x, N_\mu y)= 0$$
for all $x, y\in V\otimes \hat \Q^{\square}$. This $N_{\mu}$ is $B$-linear if and only if  $\mu(bx, y)=\mu(x, b^*y)$ for all $b\in B$ and $x,y\in V$, and in this case,  $N_\mu$   is regarded as an element of $\Lie(G) \otimes \hat \Q^{\square}$.

\end{sbpara}

\begin{sbpara}\label{cond4} 
  $\bar F_{\cH, \Sig}(S)$ is the subset of $\bar F_\cH(S)$ consisting of elements satisfying the following condition (i).
 
 \medskip

(i)  For each geometric point $s$ of $S$, let $W={\bar Y}_s\otimes \Z_{(\square)}$, and let $g: V \otimes \hat \Q^{\square}\to W \otimes \hat \Q^{\square}$ be the surjective homomorphism induced by the inverse of the level structure of $A$. 
  Then 
there exists a $\sig\in \Sig(W,g)$ such that for every $a\in \Hom((M_S/\cO^\times_S)_s,\N)$, the composition $W \times W \to  (M_S^{\gp}/\cO_S^\times)_s \otimes \Q\overset{a}\to \R$ belongs to $\sig$.

This  (i)   for each $s$ is equivalent to the following condition (i($s$)) concerning the local monodromy  at $s$ of the log abelian variety $A$. 

(i($s$)) Fix an isomorphism $V^{\square}A\cong V\otimes \hat \Q^{\square}$ which belongs to the level structure of $A$. Then there exists a $\sig\in\Sig(W, g)$ such that 
for every $\gamma\in\pi^{\log}_1(s)^+:=\Hom((M^{\gp}_S/\cO^\times_S)_s,\N)\subset \pi_1^{\log}(s)(-1)$, $\log(\gamma): V^{\square}A \to (V^{\square}A)(-1)$ (\ref{cond4-0} (2)) coincides with $N_\mu$ (\ref{cond4-0}  (3)) for some $\mu$ in $\sig\cap \Hom({\bar Y}_s\otimes{\bar Y}_s, \Q)$.
\end{sbpara}

\begin{sbpara}\label{locmod}  
Moduli functor for a single $\sigma$. 

Fix a finitely generated $\cO_{(\square)}$-module $W$ which is torsion-free as a $\Z$-module and a surjective $B\otimes \hat \Q^{\square}$-homomorphism $g: V\otimes\hat \Q^{\square} \to W \otimes\hat \Q^{\square}$ whose kernel $J$ satisfies 
$J^{\perp} \subset J$. 
  Fix also $\sig \in \Sig(W,g)$. 

Then we define $\bar F'_{\cH, \sig}(S)=$ a pair of (the same thing as in the definition of $\bar F_\cH(S)$) and (a surjective $\cO_{(\square)}$-homomorphism $g': W\to \bar Y \otimes \Z_{(\square)}$) such that the composition $V\otimes \hat \Q^{\square}\overset{g}\to W \otimes \hat \Q^{\square}\overset{g'}\to \bar Y\otimes \hat \Q^{\square}$ is induced by the inverse of the level structure and 
such that for every geometric point $s$ of $S$ and every $a\in \Hom((M_S/\cO^\times_S)_s,\N)$, the bilinear form 
$W \times W \to \bR$ induced by $g'$ and $a$ belongs to $\sig$.

We will show that $\bar F'_{\cH, \sig}$  is represented by an algebraic space   with an fs log structure (Theorem \ref{t:Main_sig}).

\end{sbpara}

\subsection{More integral moduli functors}\label{s:integral}

\begin{sbpara}\label{BU} We consider more integral structures of the objects in \ref{BV}. 
  Then we can rewrite the moduli functors in Sections \ref{s:functor} and \ref{s:fan} using the usual categories of abelian varieties without isogeny, and of log abelian varieties without isogeny.

  Let the notation and the assumptions be as in Section \ref{s:functor}.
\end{sbpara}

\begin{sbpara}\label{BU1} 

Assume that we are given a $\Z$-subalgebra $\cO$ of $\cO_{(\square)}$ such that $\cO$ is finitely generated as a $\Z$-module, 
$\cO^*=\cO$, and such that $\cO_{(\square)}= \cO\otimes \Z_{(\square)}$, and a 
 finitely generated $\cO$-submodule $L$ of $L_{(\square)}$ such that $L_{(\square)}=L\otimes \Z_{(\square)}$ and  such that $\psi(L, L)\subset \Z(1)$.
  (Now the notation and the assumptions are the same as those in \cite{Lan} 1.4.1.) 

  Then $G$ extends to an affine group scheme  $G_L$ over $\bZ$ defined as follows. For a commutative ring $R$ over $\bZ$, 
$$G_L(R)= \{(a,c)\in \text{Aut}_{R\otimes\cO}
(R\otimes L) \times R^\times\;|\; \psi(ax,ay)= c\psi(x,y) \;\text{ for all}\; x,y\}.$$

 \end{sbpara}
 
 \begin{sbrem}\label{remOL}  Such $\cO$ and $L$ such that $L\otimes \hat \Z^{\square}$ is stable under $\cH$ exist. In fact, we can easily find $L$ and then $\cO$ is obtained as 
 $\cO'\cap (\cO')^*$, where $\cO'=\{b\in \cO_{(\square)}\;|\; bL\subset L\}$.

 \end{sbrem}

\begin{sbpara}\label{FK1} In the case $\cH\subset 
G_L(\hat \Z^{\square})$ ($\cH$ is assumed to be neat),
the  moduli functors $$F_\cH : (\text{sch}/O_{F_0,(\square)})\to (\mathrm{set}), \quad  \bar F_\cH: (\fs/O_{F_0, (\square)})\to (\mathrm{set})$$
in Section \ref{s:functor} are understood  as in \ref{FK2} below. 
\end{sbpara}

\begin{sbpara}\label{FK2} 
  We have $F_\cH(S)\; (\text{resp.}\; \bar F_\cH(S))= \{(A, \iota, p, \eta)\}$.

0. $A$ is an abelian scheme (resp.\ a log abelian variety) over $S$ of dimension $g$.  (Here not like in Section \ref{s:functor}, this is considered not in the category of isogeny.)

1. $\iota$ is a homomorphism $\cO\to \text{End}(A)$ satisfying a certain condition. 

2. $p$ is a polarization whose degree is prime to $\square$ satisfying a certain condition. 

3. $\eta$ is a class of an $\cO\otimes \hat \bZ^{\square}$-isomorphism 
modulo $\cH$ on the pro-\'etale (resp.\ pro-k\'et) site of $S$ from $T^{\square}A:= \prod_{\ell\notin \square}  T_{\ell}A$ to $L \otimes \hat \Z^{\square}$ satisfying a certain condition.  

These certain conditions are similar to those in Section \ref{s:functor}, explained below.  
\end{sbpara}

\begin{sbpara}\label{cond12}
  The certain conditions in 
more integral forms are as follows. 

  The condition in \ref{FK2} 1 is that $\det_{\cO\vert\Lie_{A/S}}$ agrees with the image of $\det_{\cO \vert V_0}$, where 
$\det_{\cO\vert\Lie_{A/S}}$ and $\det_{\cO \vert V_0}$ are defined similarly as in \ref{cond1}. 

  The condition in \ref{FK2} 2 
is that for $b\in \cO$, the following diagram is commutative:
$$\begin{matrix}  A  & \overset{b^*}\to & A \\
{}^p \downarrow {\;}&& {\;}\downarrow {}^p\\
\cE xt(A, \Gmlog) & \overset{\text{by $b$}}\to & \cExt(A, \Gmlog).\end{matrix}$$

 The condition in \ref{FK2} 3 
 is that $T^{\square}A \times T^{\square}A\overset{1\times p}\to T^{\square}A \times T^{\square}\cExt(A,\Gmlog)\to \hat \Z^{\square}(1)$ is compatible with $\psi$.
\end{sbpara}

  To prove that the moduli functors are expressed in the above way, 
  we need two lemmas. 

\begin{sblem}
\label{l:tate}
  Let the notation be as above. 
  Let $A$ be a log abelian variety over $S$. 
  Let $T$ be a $\hat \Z^{\square}$-lattice in $V^{\square}A$. 
  Then there is another log abelian variety $A'$ over $S$ such that $A$ and $A'$ coincide in the
   category modulo isogenies whose degrees are prime to $\square$ and such that
 the image of the natural homomorphism $T^{\square}A' \to V^{\square}A$ coincides with $T$. 
\end{sblem}

\begin{pf}
  We may assume $nT^{\square}A \subset T \subset T^{\square}A$ for some $n$ which is prime to $\square$. 
  Let $H$ be the image of $T$ in $A[n]$. 
  It is enough to show that $A':=A/H$ is a log abelian variety. 
  Let $\bar X \times \bar Y \to \Gmlog/\Gm$ and $\cG$ be the admissible pairing and the semiabelian part of $A$, respectively. 
  Let $H_1$ be the kernel of the composition $H \to (\text{torsion of}\; A) \to {\bar Y} \otimes \Q/\Z$.  
  Let  $\cG'=\cG/H_1$. 
  Let $\bar Y'\supset \bar Y$ be the kernel of $\bar Y \otimes \Q \to \text{Coker}(H \to {\bar Y} \otimes \Q/\Z)= ({\bar Y} \otimes \Q/\Z)/(H/H_1)$.
 As in \cite{KKN2}, $\bar X$ is the character group of $\cG$ in the sense of \cite{FC} Chapter I, Theorem 2.10. 
  Let $\bar X'$ be the character group of $\cG'$. 
  
  Consider the definition \cite{KKN2} 4.1 of log abelian variety consisting of the conditions 4.1.1--4.1.3. Concerning the condition 4.1.1, the pullback of $A'$ to a geometric point $s$ of $S$ is a log abelian variety because over $s$, $A'$ corresponds to the log $1$-motif (\cite{KKN2} Section 2) $Y'\to \cG'_{\log}$ which is induced from the log $1$-motif $Y\to \cG_{\log}$ corresponding to $A$, endowed with the polarization induced by the polarization of $Y \to \cG_{\log}$. Concerning the condition 4.1.2, we have
 an exact sequence $0\to \cG' \to A' \to \cHom({\bar X'}, \Gmlog/\Gm)^{(\bar Y')}/\bar Y' \to 0$.
  The condition 4.1.3 of the separability follows from the fact that $A'$ is an $A[n]/H$-torsor over $A$ via  $n:A'\to A$. 
\end{pf}
    
\begin{sblem}
\label{l:Hom}
  Let the notation be as above. 
  Let $A$ and $A'$ be abelian schemes (resp.\ log abelian varieties) over $S$.   Then we have 

 $$\Hom(A, A') \overset{\cong}\to \{h\in \Hom(A, A') \otimes \Z_{(\square)} \; |\; h(T^{\square}A) \subset T^{\square}(A') \; \text{in}\; V^{\square}(A')\}.$$
\end{sblem}

\begin{pf} Let $h$ be an element of the set on the right hand side. Then $h=f/n$ for some $f\in \Hom(A, A')$ and for some integer $n\geq 1$ which is prime to $\square$, and we have $f(T^{\square}A)\subset n T^{\square}(A')$. The last property of $f$ shows that $f: A\to A'$ kills the $n$-torsion part $A[n]$ of $A$. Because we have an exact sequence
 $$0\to A[n] \to A \overset{n}\to A \to 0$$
 of the sheaves of abelian groups on the k\'et site, $f=ng$ for some $g\in \Hom(A, A')$. Hence $h\in \Hom(A, A')$. 
\end{pf}

\begin{sbrem}
  A related result is as follows: Assume that, for each $j=1,2$, log abelian varieties $A_j$ and $A_j'$ coincide in the category modulo isogenies whose degrees are prime to $\square$.
  Let $m, n$ be  integers which are prime to $\square$ such that  $nT^{\square}A_j\subset m T^{\square}A_j'\subset T^{\square}A_j$. 
  Identify $m T^{\square}A'_j/nT^{\square}A_j$ with a subgroup $H_j$ of $A_j[n]$. 
  Then we have $\Hom(A_1', A_2')= \{h\in \Hom(A_1, A_2)\;|\; h(H_1)\subset H_2\}$. 
  This follows from Lemma \ref{l:Hom}.  It also follows from the fact that $A_j'$ is identified with $A_j/H_j$. 
\end{sbrem}

\begin{sbpara}
  For a given $(A',\iota, p, \eta)\in F_\cH(S)$ (resp.\ $\bar F_\cH(S)$) in Section \ref{s:functor}, 
  we see that $A'$ gives an abelian scheme (resp.\ a log abelian variety) $A$ such that the inverse image of $L\otimes \hat \Z^{\square}$ in $V^{\square}A'$ under $\eta$ is $T^{\square}A$ (Lemma \ref{l:tate}), and that,   
by Lemma \ref{l:Hom}, $F_\cH(S)$ (resp.\ $\bar F_\cH(S)$) in Section \ref{s:functor} is described as in the above way in terms of $A$ not using isogeny. 
\end{sbpara}

\begin{sbpara}
Accordingly, for $\Sig$ in Section \ref{s:fan}, the definitions of $\bar F_{\cH, \Sig}$ and $\bar F'_{\cH, \sig}$ ($\cH\subset G_L(\hat \Z^{\square})$, $\sig\in \Sig$) can be described by using the above presentation of  $\bar F_\cH$.  

\end{sbpara}

\begin{sbrem} As above, more integral presentations of $F_\cH$ and $\bar F_\cH$ without working modulo isogeny are given in the case $\cH\subset G_L(\hat \Z^{\square})$. 

However by Remark \ref{remOL}, 
for any open compact subgroup $\cH$ of $G(\hat \Q^{\square})$, by taking $\cO$ and $L$ such that $\cH\subset G_L(\hat \Z^{\square})$, we can express $F_\cH$ and $\bar F_\cH$ in the more integral ways without using isogeny. Thus for the proofs of the main results of this paper, we can use the more integral formulations of the moduli functors given in this section \ref{s:integral}. This is important because the representability of $\cH om(A, A')$ (explained in Section \ref{s:hom}) works nicely for the more integral formulations. 
\end{sbrem}

\section{The main theorems and proofs}
\label{s:main}

\subsection{Log algebraic spaces (review)}

\begin{sbpara}
  Before stating the main results, we recall the definitions of log algebraic space in the first sense and in the second sense from \cite{KKN4} 10.1. 

  Let $S$ be an fs log scheme. 

  A {\it log algebraic space over $S$ in the first sense} is an algebraic space (in the usual sense) endowed with an fs log structure. 
  We identify it with the sheaf on $(\fs/S)_{\et}$ represented by it. 

  A sheaf $F$ on $(\fs/S)_{\et}$ is a {\it log algebraic space in the second sense} over $S$ if there is a surjective morphism $F' \to F$ of sheaves from a log algebraic space $F'$ in the first sense such that for any morphism $T \to F$ from an fs log scheme $T$, the base change $F'\times_FT \to T$ is a log \'etale log algebraic space over $T$ in the first sense.
\end{sbpara}

\begin{sbpara}
  A log algebraic space in the first sense is a log algebraic space in the second sense (\cite{KKN4} Proposition 10.2). 
  A log abelian variety is a log algebraic space in the second sense (\cite{KKN4} Theorem 10.4 (1)), though it is not necessarily so in the first sense. 
\end{sbpara}

\subsection{Main results}\label{s:Main}

We state the main results in this paper. 
  Let the notation and the assumptions be as in Section \ref{s:mod}.

\begin{sbthm}
\label{t:Main} 
The moduli functor $\bar F_\cH$ (without $\Sig$) is represented by a proper and log smooth log algebraic space  in the second sense. 
\end{sbthm}

\begin{sbthm}
\label{t:Main_Sig}
The moduli functor $\bar F_{\cH,\Sig}$ (with $\Sig$) is represented by a proper and log smooth log algebraic space  in the first sense.
\end{sbthm}

\begin{sbthm}
\label{t:Main_sig} 
$\bar F'_{\cH, \sig}$ is a log smooth log algebraic space in the first sense. \end{sbthm}

\subsection{Homomorphisms of log abelian varieties}
\label{s:hom}

In this section, we prove the following proposition.

\begin{sbprop}\label{hom0} 
  Let $S$ be an fs log scheme and let $A$ and $A'$ be principally polarizable log abelian varieties over $S$. Then $\cHom(A, A')$ is represented by an algebraic space over $S$ endowed with the inverse image of the log structure of $S$.
\end{sbprop}

  First, our former works give some good properties of this functor. 

\begin{sbpara}\label{lfp} By \cite{KKN4} Proposition 9.2 (1), $\cHom(A, A')$ is locally of finite presentation. 
\end{sbpara}

\begin{sbpara}\label{hom2}
GAGF for $\cHom(A,A')$ holds.  See \cite{KKN6} Theorem 6.1 and Remark 6.1.1.
\end{sbpara}

\begin{sbpara}\label{hom1} We prove Proposition \ref{hom0} in the case where $A$ and $A'$ are with constant degenerations. 

  Let $Y, \cG, \cT, \cB$ be the lattice, the semiabelian part, the torus part, and the abelian part of $A$ respectively, and let $Y', \cG', \cT', \cB'$ be those for $A'$ respectively. 
  Recall that we have an exact sequence 
$0 \to \cT \to \cG \to \cB \to 0$ and the presentation $A=\cG_{\log}^{(Y)}/Y$. 
  Note that $Y$ is a locally constant sheaf because $A$ is with constant degeneration. 

  It is enough to prove that the canonical homomorphism
\begin{equation*}\tag{$*$}
\cHom(A, A')   \to   \cHom(Y, Y') \times \cHom(\cT, \cT') \times \cHom(\cB, \cB')
\end{equation*}
is represented by strict closed immersions. 

First,  $\cHom(\cG, \cG')$ is the fiber product of  $\cHom(\cT, \cT') \times \cHom(\cB, \cB') \to E \times E \leftarrow E$, where  $E=\cExt(\cB, \cT')$ and the last arrow is the diagonal map. Note that $E$ is an abelian scheme over $S$. 

Next,  $\cHom(A, A')$ is the fiber product of  $\cHom(Y, Y') \times \cHom(\cG, \cG') \to L \times L \leftarrow L$, where  $ L= \cHom(Y, (\cG')_{\log}^{(Y')})$ and the last arrow is the diagonal map.  This arrow $L \to L \times L$ is  represented by strict closed immersions.

  Therefore the homomorphism $(*)$ is represented by strict closed immersions. 
\end{sbpara} 

\noindent
{\it Remark.} 
  The fact that in the notation and the assumptions in Proposition \ref{hom0}, for a noetherian $S$, $\Hom(A, A')$ (not $\cHom(A, A')$) is a finitely generated abelian group can be deduced from the above argument by considering the log constant loci of $S$.  

\begin{sbpara}
\label{strict}
  Let $U'\to U$ be a morphism of fs log schemes over $S$ whose underlying morphism of schemes is an isomorphism. 
  Then the homomorphism $\cHom(A,A')(U) \to \cHom(A,A')(U')$ is bijective.
  To see it, by \ref{lfp} and \ref{hom2}, we may assume that $S$ is the spectrum of an Artin ring (cf.\ \cite{KKN7} Remark 6.5). 
  In this case, $A$ and $A'$ are with constant degeneration. 
  Then a homomorphism between them is determined by the pair consisting of a homomorphism of the lattices and a homomorphism of the semiabelian parts (cf.\ \cite{KKN2} Proposition 2.5).
  Hence we see that the above homomorphism is bijective. 
\end{sbpara}

\begin{sbpara}\label{hom12}
By \ref{strict}, 
we have $\cHom(A,A')(U^{\mathrm{str}})=\cHom(A,A')(U)$. 
  Here $U$ is an fs log scheme over $S$ and $U^{\mathrm{str}}$ is the fs log scheme over $S$ whose underlying scheme is that of $U$ and whose log structure  is  the inverse image of the log structure of $S$.
\end{sbpara}

\begin{sbpara}
\label{def_theory}
  Below, we consider 
the functor associating to a scheme $U$ over the underlying scheme of $S$ the set 
$\Hom_U(A,A')$, where $U$ is endowed with the inverse image of the log structure of $S$.
  We denote this functor by $H$. 

  By \ref{hom12}, to prove Proposition \ref{hom0}, it suffices to show that $H$ is represented by an algebraic space. 
  Since we want to use Artin's criterion, we have to know the deformation of $H$. 

\end{sbpara}

By the theory of universal additive extensions (\cite{KKN7} Section 3), 
we have the following. 

\begin{sbpara}\label{hom3} Let $U_0$ be an fs log scheme over $S$ and let $U_1$ be a nilpotent thickening of $U_0$ over $S$ such that the kernel $I$ of $\cO_{U_1}\to \cO_{U_0}$ satisfies $I^2=0$.
Then $\Hom(A_{U_1}, A'_{U_1})$ is identified with the subgroup of $\Hom(A_{U_0}, A'_{U_0})$ consisting of $h$ such that the homomorphism 
$\Lie(E(A))_{U_1}\to \Lie(E(A'))_{U_1} $ induced by $h$ respects the filtrations. 
  Here $E(A)$ (resp.\ $E(A')$) denotes the universal additive extension and 
the filtration on $\Lie(E(A))_{U_1}$ (resp.\ $\Lie(E(A'))_{U_1}$) is the two-step one determined by the kernel of the map 
$\Lie(E(A))_{U_1} \to  \Lie(A_{U_1})$ (resp.\ 
$\Lie(E(A'))_{U_1} \to \Lie(A'_{U_1})$).
\end{sbpara}

\begin{sbpara}\label{hom31} By \ref{hom3}, we see that  $\cHom(A, A')$ is unramified.

\end{sbpara}

\begin{sbpara}\label{obs1} 
  Let the notation be as in \ref{hom3}. 
We show that we have an exact sequence 
$$0\to \Hom(A_{U_1}, A'_{U_1})\to \Hom(A_{U_0}, A'_{U_0}) \to \Hom(\coLie(A^*_{U_0}), \Lie(A'_{U_0})\otimes_{\cO_{U_0}} I),$$
where $(\cdot)^*$ denotes the dual log abelian variety (cf.\ \cite{KKN7} 3.3).

Let $h\in \Hom(A_{U_0}, A'_{U_0})$. Then $h$ induces 
$l: \Lie(E(A))_{U_1}\to \Lie(E(A'))_{U_1}$. 
  Since we have an exact sequence $$0\to \coLie(A^*_{U_1}) \to \Lie(E(A))_{U_1} \to \Lie(A_{U_1})\to 0$$ 
and a similar exact sequence for $A'$,  $l$ respects the filtrations if and only if the map $\delta: \coLie(A^*_{U_1}) \to \Lie(A'_{U_1})$ induced by $l$ is the zero map.
Since the map $\Lie(E(A))_{U_0}\to \Lie(E(A'))_{U_0}$ induced by $l$ respects the filtrations, the image of $\delta$ is contained in  $I\cdot \Lie(A'_{U_1}) = \Lie(A'_{U_0})\otimes_{\cO_{U_0}} I$ and $\delta$ factors as $\coLie(A^*_{U_1})\to \coLie(A^*_{U_0}) \to \Lie(A'_{U_0}) \otimes_{\cO_{U_0}} I$. 
  This gives the desired exact sequence. 
\end{sbpara}

\begin{sbpara}\label{fibp} 
  Let $P\to R$ be a surjective ring homomorphism over $S$ whose kernel is a nilpotent ideal, and let $Q\to R$ be a ring homomorphism over $S$. 
  We claim that for $H$ in \ref{def_theory}, 
the map $H(P\times_R Q)\to H(P) \times_{H(R)} H(Q)$ is bijective. 

  Since the homomorphism $P\to R$ can be decomposed into 
$P=P_n \to P_{n-1}\to \dots \to P_0=R$ such that $P_i \to P_{i-1}$ are surjective with square zero kernels for all $i$, we may assume 
that the kernel $I$ of $P\to R$ satisfies $I^2=0$. 
  In this case, the kernel of $P\times_R Q \to Q$ is also $I$. 
  Hence by \ref{obs1}, we have a commutative diagram of exact sequences
$$\begin{matrix}
0& \to & H(P \times_R Q) & \to & H(Q) & \to & \Hom_Q(\coLie(A^*_Q), \Lie(A_Q)\otimes_Q I) \\
&&\downarrow &&\downarrow &&\Vert\\
0& \to & H(P) & \to & H(R) & \to & \Hom_R(\coLie(A^*_R), \Lie(A_R)\otimes_R I).
\end{matrix}$$
This proves our claim. 
\end{sbpara}

\begin{sbpara} From these results, we can deduce Proposition \ref{hom0} by using Artin's criterion. 

We can use the original theorem in \cite{Artin:criterion}, 
or we can use the theorem of Hall and Rydh in \cite{HR}.
\end{sbpara}

\begin{sbpara} We can check Artin's original conditions in 
\cite{Artin:criterion} Theorem 5.3. 
  We may assume that the base $S$ is finitely generated over $\bZ$. 

  Below, the notation follows loc.\ cit.
  As a deformation theory for $H$, we take the standard one: $D(A_0,M,\xi_0):=H_{\xi_0}(A_0[M])$, which vanishes by \ref{hom31} so that the conditions  [4${}'$](a)--(c) are trivially satisfied.

  Since $A'$ is an fppf sheaf (cf.\ \cite{KKN7} 6.3), the condition [0${}'$] is satisfied. 

  By \ref{lfp} and \ref{hom2}, the conditions [1${}'$] and [2${}'$] are satisfied, respectively. 

  By \ref{strict}, [3${}'$](a) is reduced to \cite{KKN6} Theorem 3.4 (cf.\ \cite{KKN7} 6.8).

  To see that the conditions [3${}'$](b) and [5${}'$](c) are satisfied, we may assume that $A$ and $A'$ are with constant degeneration (cf.\ \cite{KKN7} 6.10) so that they reduce to \ref{hom1}. 

 Finally, by \ref{fibp}, the conditions [5${}'$](a) and [5${}'$](b) are satisfied.
\end{sbpara}

\begin{sbpara} We can also check the conditions of  Hall and Rydh in \cite{HR} for the representability. We can use the condition 8.3 in \cite{HR} there because the DVR-homogeneity is satisfied by  \ref{fibp}. 
  This condition 8.3 in \cite{HR} is satisfied because by \ref{hom3}, we have a canonical injective homomorphism
$\text{Obs}(U,I) \to \Hom_{\cO_U}(\coLie(A^*_U), \Lie(A'_U)\otimes_{\cO_U} I)$. 
\end{sbpara}

\subsection{Proof of the main theorem}
\label{s:pf1}

  We prove the main theorems (Theorems \ref{t:Main}, \ref{t:Main_Sig}, and \ref{t:Main_sig}).

\begin{sbpara} 
\label{pf0}
  We prove the part of Theorem \ref{t:Main} 
that 
the moduli functor $\bar F_\cH$ (without fixing $\Sig$) is a log algebraic space in the second sense.

We may assume that $\cH$ is sufficiently small. 
  This is because if $\bar F_{\cH'}$ is a log algebraic space in the second sense for a smaller $\cH'\subset \cH$, then $\bar F_\cH$ is also because $\bar F_{\cH'}$ is a k\'et finite covering of it.

  First we prove the non-coefficient case $B=\Q$. 
  This is essentially contained in \cite{KKN7} Theorem 1.6. 
  In fact, in this case, take an $L$ as in \ref{s:integral} such that $\psi$ induces a perfect pairing on $L$ (cf.\ \cite{KKN1} Example 4.2.3). 
  Then for an integer $n$ which is prime to $\square$, 
the moduli functor $\bar F_\cH$ for $\cH=\text{Ker}(G_L(\hat \Z^{\square})\to G_L(\Z/n\Z))$ is a closed open subfunctor of $F_{g,n}$ in \cite{KKN7} (cf.\ \cite{KKN7} 9.1).
  Since $F_{g,n}$ is a log algebraic space in the second sense by the above theorem in \cite{KKN7}, $\bar F_\cH$ is also a log algebraic space in the second sense, and it holds for any $\cH$ by the previous paragraph. 

  The general case follows from the non-coefficient case and the representability of $\cH om(A, A')$ (Proposition \ref{hom0}) just as in the non-log case treated, for example, in \cite{GN} Theorem 3.4.3. 
  In fact, take an $L$ such that $L \otimes \hat \Z^{\square}$ is $\cH$-stable, and take an open neat subgroup $\cH'$ of $GSp_L(\hat \Z^{\square})$ such that $\cH= G_L(\hat \Z^{\square})\cap \cH'$. 
  We consider the moduli space $F$ with no coefficients with respect to these $L$ and $\cH'$. 
  Since we already know that $F$ is a log algebraic space in the second sense, there are an fs log scheme $X$ and a surjective morphism $X \to F$ which is relatively represented by log \'etale log algebraic spaces in the first sense. 
  Consider the projection $\bar F_\cH \to F$ (forgetting the coefficients) and the fiber product 
$\bar F_\cH\times_FX$. 
  Applying Proposition \ref{hom0} to the pullback $A$ to $X$ of the universal family on $F$, we see that to add 
the coefficient $\cO \to \End(A)$ is representable and that 
$\bar F_\cH\times_FX$ is a log algebraic space in the first sense. 
  Here, to give an $\cH$-level structure gives a closed open subfunctor, because $G_L(\hat \Z^{\square})/\cH$  is a  subset of the finite set $GSp_L(\hat \Z^{\square})/\cH'$.  
  (Also we see that the compatibility conditions (\ref{cond12}) are harmless as in the first part of the proof of \cite{GN} Theorem 3.4.3 
by taking 
\cite{KKN4} Proposition 18.1 into account for the condition 3, and 
by observing that the condition 2 is equivalent to $p \circ b^* = b \circ p: A \to A^*$.)
  Thus $\bar F_\cH$ is a log algebraic space in the second sense. 
\end{sbpara}

\begin{sbpara}\label{mainpf1}   
  In the same way as in \ref{pf0}, 
the part of Theorem \ref{t:Main_sig}
that $\bar F'_{\cH, \sig}$ is a log algebraic space in the first sense 
is reduced to the non-coefficient case (\cite{KKN7} 7.2) by the representability of $\cHom$ (Proposition \ref{hom0}). 
\end{sbpara}

\begin{sbpara}
\label{mainpf2}
The part of Theorem \ref{t:Main_Sig} 
that $\bar F_{\cH, \Sig}$ is a log algebraic space in the first sense 
is reduced to the part of Theorem \ref{t:Main_sig} proved in \ref{mainpf1} by using the strict \'etale morphisms $\bar F'_{\cH, \sig}\to \bar F_{\cH, \Sig}$ for all $\sig$ in all $\Sig(W,g)$. 

  Note that the part of Theorem \ref{t:Main} except the properness and the log smoothness is also deduced from \ref{mainpf1} as in the non-coefficient case (cf.\ \cite{KKN7} 7.4). 

  On the other hand, in Section \ref{s:pf}, we will give an alternative proof of the above part of Theorem \ref{t:Main_Sig} without using \ref{mainpf1}. 
\end{sbpara}

\begin{sbpara} 
\label{smooth}
Proof of the log smoothness of $\bar F_{\cH, \Sig}$.

It is enough to show that 
an $(A,\iota,p,\eta) \in \bar F_{\cH,\Sig}(S)$ as in \ref{FK22} 
can be lifted to a nilpotent thickening of $S$. 
We may assume that (up to isogenies as in \ref{FK22}) $(A,p)$ comes from a principally polarized log abelian variety. 
In fact, take any representative $(A,p)$ in the category without isogenies. 
Then, forgetting the coefficients, we can take a lattice $T$ of $V^{\square}A$ on which $\psi$ and $\eta$ induce a perfect pairing (cf.\ the third paragraph of \ref{pf0}). 
  By Lemma \ref{l:tate}, we find the desired representative as the one corresponding to $T$. 
\end{sbpara}

  Hence the log smoothness is essentially reduced to the following. 
  (See the complement after the lemma for the determinant condition.)

\begin{sblem}
  Let $S\to S'$ be a strict nilpotent thickening of affine fs log schemes of finite type over $\bZ$ such that the ideal $I$ of $\cO_{S'}$ defining $S$ satisfies $I^2=0$. 
  Let $A$ be a principally polarized log abelian variety over $S$ with 
a homomorphism $\cO_{(\square)} \to \End(A)\otimes \bZ_{(\square)}$. 
  Then, Zariski locally on $S$, there is a principally polarized log abelian variety $A'$ over $S'$ with 
a homomorphism $\cO_{(\square)} \to \End(A')\otimes \bZ_{(\square)}$, whose restriction to $S$ coincides with the given one. 
\end{sblem}

\begin{pf}
  In this proof, we call 
a homomorphism $\cO_{(\square)} \to \End(A)\otimes \bZ_{(\square)}$
an action of $\cO_{(\square)}$. 

  Locally on $S$, take a lift $A''$ of $A$ to $S'$ (\cite{KKN7} Proposition 4.2) which may not have an action of $\cO_{(\square)}$. 
  Then the universal additive extension $E(A'') $ of $A''$ does not depend on the choice of $A''$ (\cite{KKN7} Proposition 3.7) and the action of $\cO_{(\square)}$ on $A$ induces a  canonical action of $\cO_{(\square)}$ on $E(A'')$. The kernel of $E(A'')\to A'' $ is $\coLie((A'')^*)$. The action of $\cO_{(\square)}$ on $E(A'')$ keeps the kernel $K$ of $E(A'')\to \Lie(A)$. We have an exact sequence $0\to \coLie((A'')^*)\to K \to \Lie(A) \otimes_{\cO_{S}} I \to 0$. For $b\in \cO_{(\square)}$ and $x\in \coLie((A'')^*)$, let $\delta(b)x$ be the image of $bx\in K$ in $\Lie(A) \otimes_{\cO_S} I$. This defines a homomorphism $\delta: \cO_{(\square)} \to H:=\cH om(\coLie(A^*), \Lie(A))\otimes_{\cO_S} I$. This $\delta$ is a derivation, that is, $\delta(bb')= b\delta(b')+ \delta(b)b'$ for all $b, b'\in \cO_{(\square)}$. Since $\cO_{(\square)}$ is an Azumaya algebra over its center, which is finite \'etale over $O_{F_0, (\square)}$, this derivation $\delta$ is an inner-derivation, that is, there is an $h\in H$ such that $\delta(b)= bh-hb$ for all $b\in \cO_{(\square)}$. Let the lift $A'$ of $A$ to $S'$ be the twist of $A''$ by this element $h$ of $H$ (\cite{KKN7} Proposition 4.1). Then the action of $\cO_{(\square)}$ on $E(A'')$ keeps the kernel of $E(A'')\to A'$ and hence $\cO_{(\square)}$ acts on $A'$. 
\end{pf}

Thus we have a nice lift $A'$ of $A$ to $S'$. 
  The moduli condition on the determinant of the action of $\cO_{(\square)}$ on $\Lie(A')$ is automatically satisfied.
  In fact, \'etale locally on $\Spec(\cO_{F_0, (\square)})$, $\cO_{(\square)}$ becomes a finite product of matrix algebras, and by using the Morita equivalence between modules over a ring and modules over the matrix ring over the ring, the condition concerned for $A'$ is easily proved.

\begin{sbpara}
  Similarly, 
  we 
see the log smoothness of $\bar F'_{\cH,\sig}$, 
which completes the proof of Theorem \ref{t:Main_sig}, and from which we can see the log smoothness of $\bar F_\cH$ as in \cite{KKN7} 7.4.
\end{sbpara}

\begin{sbpara}
\label{proper}
  We prove the properness of $\bar F_{\cH}$ and $\bar F_{\cH, \Sig}$, which completes the proofs of Theorems \ref{t:Main} and \ref{t:Main_Sig}.
  By the non-coefficient case in \cite{KKN7} Theorems 1.6 and 1.7, it is reduced to the properness of $\cHom(A,A')$ in Proposition \ref{hom0}. 
  The latter is shown as in \cite{KKN7} Section 8 by the valuative criterion and the equivalence between the category of log abelian varieties over a complete discrete valuation field and that of ones over its valuation ring.
\end{sbpara}

\section{Second proof}
\label{s:artin}

  In this section, we give the second proof of the main results in this paper by using log Artin's criterion together with other applications of this criterion.

\subsection{Log Artin's criterion}\label{s:lAr}

We improve  \cite{KKN4} Proposition 7.2 on log Artin's criterion. 

\begin{sbthm}\label{logAr} 
  Let $S$ be an fs log scheme 
whose underlying scheme is of finite type over a field or over an excellent Dedekind domain. 
  For a sheaf $F$ on $(\fs/S)_{\et}$, the following two conditions are equivalent.

{\rm (i)} $F$ is represented by an algebraic space with an fs log structure over $S$ (that is, $F$ is a log algebraic space in the first sense). 

{\rm (ii)} $F$ satisfies the conditions $7.2.1$--$7.2.4$ in 
{\rm\cite{KKN4}} Proposition $7.2$, that is, the quasi-separatedness, the finiteness, the effective pro-representability, and the openness of versality. 
\end{sbthm}

\noindent {\it Remark.} 
  M.\ C.\ Olsson proves another version of log Artin's criterion (\cite{Olsson} Section 3.5). 

\begin{pf} In \cite{KKN4} Proposition 7.2, it was shown that (i) implies (ii) and that (ii) plus a monomorphism condition in the latter half of (1) in \cite{KKN4} 6.2 implies (i).  We show that actually, (ii) implies this monomorphism condition. 

Let $S$ be an fs log scheme (resp.\ a scheme). 
For a functor $F$ on $(\fs/S)$ (resp.\ $(\text{sch}/S)$), let $[[F]]$ be the restriction of $F$ to 
the subcategory consisting of objects whose underlying schemes (resp.\ which) are the spectrums of Artin rings.   
  If $F$ is represented by a locally noetherian  fs log scheme $X$ (resp.\ locally noetherian scheme $X$) over $S$, $[[F]]$ is pro-represented by the disjoint union of $\text{Spf}(\cO_{X,x})$ for $x\in X$. (In the log situation, this Spf is endowed with the inverse image of the log structure of $X$.)
In the following, for an fs log scheme $X$ (resp.\ an fs log formal  scheme) $X$, $X^{\circ}$ denotes the underlying scheme (resp.\ underlying formal scheme). We have 

\medskip

(\ref{logAr}.1) For a locally noetherian fs log scheme $X$, we have $[[X^{\circ}]]= [[X]]^{\circ}$.

\medskip

Let $F$ be a functor satisfying the condition (ii). 
  By \cite{KKN4} Proposition 7.2, we have an fs log scheme $F'$ and a strict \'etale surjective morphism $F'\to F$. 

The monomorphism condition in problem is that 

\medskip

(\ref{logAr}.2) \quad $(F' \times_F F')^{\circ} \to (F' \times_{S^{\circ}} F')^{\circ}$

\medskip

\noindent 
is a monomorphism in the category of schemes. 
This is reduced to that the morphism

\medskip

(\ref{logAr}.3) \quad $[[(F' \times_F F')^{\circ}]] \to [[(F' \times_{S^{\circ}} F')^{\circ}]]$

\medskip

\noindent 
is a monomorphism in the category of formal schemes. (We can reduce \ref{logAr}.2 to \ref{logAr}.3 by the fact that any noetherian ring can be embedded into the inverse limit of Artin rings.) 

We rewrite the left hand side and the right hand side of \ref{logAr}.3.

First, the left hand side is rewritten as 
$$[[(F' \times_F F')^{\circ}]]=[[F' \times_F F']]^{\circ}= ([[F']] \times_{[[F]]} [[F']])^{\circ}= [[F']]^{\circ} \times_{[[F]]^{\circ}} [[F']]^{\circ}.$$
Here the first equality follows from \ref{logAr}.1, the second equality is clear, and the third equality follows from the fact that $[[F]]$ is representable (by the pro-representability in (ii)) and the fact that $[[F']]\to [[F]]$ is strict. 

The right hand side is rewritten as 
$$[[(F' \times_{S^{\circ}} F')^{\circ}]] = [[(F')^{\circ} \times_{S^{\circ}} (F')^{\circ}]]= [[(F')^{\circ}]] \times_{[[S^{\circ}]]} [[(F')^{\circ}]] = [[F']]^{\circ} \times_{[[S]]^{\circ}} [[F']]^{\circ}.$$
Here the first two equalities are evident and the last equality follows from \ref{logAr}.1. 

Hence we are reduced to that the morphism 
$$[[F']]^{\circ} \times_{[[F]]^{\circ}} [[F']]^{\circ}\to  [[F']]^{\circ} \times_{[[S]]^{\circ}} [[F']]^{\circ}$$
is a monomorphism, which is evident. 
\end{pf} 

As a consequence of the above improvement, \cite{KKN4} 6.2 is partially improved as in the following. 

\begin{sbcor}
\label{c:Ar}
  Let the notation and the assumptions be as in Theorem $\ref{logAr}$. 
  Then the conditions {\rm (i)} and  {\rm (ii)} are also equivalent to the following condition. 

$(*)$ The diagonal morphism $F \to F \times F$ is represented by morphisms of finite type, and there is a morphism $U \to F$ from an fs log scheme which is represented by surjective and strict \'etale morphisms.
\end{sbcor}

In fact, together with \cite{KKN4} 6.2, for an $F$ as above, we can give the representing log algebraic space in the first sense by $(J, M_J)$, where $J$ is defined by the underlying schemes of $U$ and $U \times_FU$, and $M_J$ is determined by the log structure of $U$. 

\subsection{Simplification of the proof of \cite{KKN7}}
\label{s:cpl3}

  We explain that the improvement of log Artin's criterion (Theorem \ref{logAr}) much  simplifies the proofs of the main results in \cite{KKN7}. 

  The hardest part of \cite{KKN7} is to show that $F'_{\sig}$ in the notation in \cite{KKN7} 7.2 is a log smooth log algebraic space in the first sense. 
  In \cite{KKN7}, we prove this by the representability of $F'_{\cal T}$ (in the notation in \cite{KKN7}). 
  But, in virtue of the improved theorem, we prove the representability of $F'_{\sig}$ directly by checking the four conditions in Theorem \ref{logAr} for $F'_{\sig}$.

\begin{sbpara}
\label{qsep}
  First, the quasi-separatedness of $F'_{\sig}$ is reduced to Proposition \ref{hom0}. 
  In fact, considering $\{(f,g)\,|\,fg=1, gf=1\} \subset \cHom(A,A') \times \cHom(A',A)$ in the notation in Proposition \ref{hom0}, we see that the diagonal for $F'_{\sig}$ is represented by a strict log algebraic space in the first sense.
  Since it is strict, the underlying morphisms of schemes are monomorphisms, hence, it is represented by schemes. 
  Further, we can check the condition [$3'$](b) in \cite{Artin:criterion} Theorem (5.3) for the underlying functor of the diagonal by reducing to the case of constant degeneration (cf.\ \cite{KKN7} 6.10).
  Hence, the diagonal is quasi-compact. 
\end{sbpara}

\begin{sbpara}
  The condition of being locally of finite presentation is essentially proved in \cite{KKN7} 6.6, where we saw that $F'$ (in the notation there) is locally of finite presentation. 
  Though we have to take care of the monodromy condition, it is easy.
\end{sbpara}

\begin{sbpara}
  The pro-representability follows from the formal moduli theory and GAGF for log abelian varieties in \cite{KKN6}.
\end{sbpara}

\begin{sbpara}
The openness of versality is proved as follows. Let $U\to F$, $u\in U$ be as in the condition. The condition of being formally strict \'etale at $u$ is equivalent to the condition of being formally strict at $u$ and formally log \'etale at $u$. 
  Let $A$ be the log abelian variety on $U$. 
  Let the notation be as in \cite{KKN7}. 

The strictness condition is restated as that the map $\bar Y_{\bar u} \otimes \bar X_{\bar u}\otimes \Q \to (M^{\gp}_U/\cO^\times_U)_{\bar u}\otimes \Q$ factors through the surjection $\bar Y_{\bar u}\otimes \bar X_{\bar u}\otimes \Q \to \Hom(\sig_{\bN}, \Q)$ and induces an isomorphism $\Hom(\tau, \N) \overset{\cong}\to (M_U/\cO^\times_U)_{\bar u}$ for some face $\tau$ of $\sig_{\bN}$, 
where $\sig_{\bN}=\{x \in \sig\,|\,x(W\times W)\subset \bZ\}$. The formal log \'etaleness condition is that $U$ is log smooth over $S$ and the map $\text{Sym}^2(\coLie(A)) \to \Omega^{1,\log}_{U/S}$ is an isomorphism at $u$. Both properties  extend to an open neighborhood of $u$.

\end{sbpara}

\begin{sbrem}
  Since the proof of the fact that $\overset \circ X$ is representable in \cite{KKN7} 7.2 (the notation is as in there) is rather rough, we give a complement. 
  In there, we checked the openness of versality of $\overset \circ X$.  
  Before it, we should see that the projection $\overset \circ H_{\sig} \to \overset \circ X$ is a torus torsor, which is seen formally. 
  Then the openness of versality of $\overset \circ X$ is deduced from that of $\overset \circ H_{\sig}.$
  The quasi-separatedness of $\overset \circ X$ is also seen formally or by the use of the above torsor property. 
\end{sbrem}

\subsection{Second proof of the main theorem}
\label{s:pf}

  Here we prove the part of Theorem \ref{t:Main_Sig} 
  except the properness and the log smoothness, which was proved in \ref{mainpf2}, by another method.  
  We do not use $\bar F'_{\cH,\sig}$ in Theorem \ref{t:Main_sig} here but 
we use log Artin's criterion (Theorem \ref{logAr}) and the part of Theorem \ref{t:Main} except the properness and the log smoothness.

\begin{sbpara}
\label{second}
  We prove that $\bar F_{\cH,\Sig}$ is a log algebraic space in the first sense. 
  We already have that $\bar F_{\cH, \Sig}$ is a log algebraic space in the second sense by the part of Theorem \ref{t:Main} proved in \ref{pf0} since $\bar F_{\cH, \Sig}\to \bar F_\cH$ is log \'etale.
\end{sbpara}

\begin{sbpara}
  By log Artin's criterion, it is enough to check the four conditions \cite{KKN4} 7.2.1--7.2.4 for $\bar F_{\cH,\Sig}$. 

  First, \cite{KKN4} 7.2.1 is reduced to Proposition \ref{hom0} similarly as in \ref{qsep}.
\end{sbpara}

\begin{sbpara}
\label{pf722}
  Next, \cite{KKN4} 7.2.2 is a consequence of \ref{second}. 
  In fact, any log algebraic space $F$ in the second sense always satisfies the condition \cite{KKN4} 7.2.2 (that is, it is locally of finite presentation). 
  This is seen as follows. 
  It can be shown that there is a surjective morphism $U\to F$ of sheaves from a log algebraic space $U$ in the first sense such that $U\times_F U$ is a log algebraic space in the first sense. This implies that \cite{KKN4} 7.2.2 for $F$ is satisfied.
\end{sbpara}

\begin{sbpara}
\label{local_moduli}
 We prove \cite{KKN4} 7.2.3 (pro-representability).
We fix an fs log point $s$ and a log abelian variety $A_0$ over $s$ with PEL structure. We consider  fs log schemes $S$ with a strict closed immersion $s\to S$ such that $S= \Spec(R)$ for a complete local ring $R$. We are discussing the moduli space (called the local moduli space) of the liftings $A$ of $A_0$ to $S$. 

The local moduli space $X$ without coefficient ring is already constructed in \cite{KKN6} Section 2. 
  To lift the coefficient, the question is whether $\cO \to \text{End}(A_0)$ induces $\cO \to \text{End}(A)\subset \text{End}(A_0)$. This is a closed condition. Hence our local moduli space is obtained as a closed subspace of $X$. 
\end{sbpara}

\begin{sbpara}\label{pf724} We prove \cite{KKN4} 7.2.4 (the openness of versality). 

Let $U\to F:=\bar F_{\cH,\Sig}$ be a morphism from an fs log scheme of finite type over $O_{F_0,(\square)}$ and $u\in U$ a point. 
  Assume that $U \to F$ is strict formally \'etale at $u$. 
  We prove that it is strict formally \'etale at any point of some open neighborhood of $u$. 
  We claim that we may assume that $U \to F$ is log \'etale. 
  In fact, since $F$ is a log algebraic space in the second sense locally of finite type (\ref{pf722}), we have a log \'etale morphism $Z \to F$ from an fs log scheme locally of finite type over $O_{F_0,(\square)}$ 
which is surjective as a morphism of sheaves for the (classical) \'etale topology.
  Then \'etale locally, the morphism $U\to F$ factors as $U\to Z \to F$. 
  Let $z$ be the image of $u$ in $Z$. 
  Then $Z\to F$ is strict formally \'etale at $z$. 
  Hence it is enough to prove that $Z\to F$ is strict on 
some open neighborhood of $z$. 
  The claim follows. 
  
  Now the condition of being strict at $u$ is equivalent to the condition that 
for some $\sig \in \Sig(W,g)$ as in (i($\bar u$)) in \ref{cond4}, the map $\pi^{\log}_1(\bar u)^+ \to \sig \cap\Hom(Y_{\bar u} \times Y_{\bar u}, \Q)\;;\;\gamma\mapsto \mu$,   $\log(\gamma)=N_\mu$, induces an isomorphism $\pi^{\log}_1(\bar u)^+\overset{\cong}\to \sig(\cH)$, where 
$$\sig(\cH)= \{\mu \in \sig\cap \Hom(Y_{\bar u} \times Y_{\bar u}, \Z_{(\square)})\;|\; \exp(N_\mu) \in \cH \;\text{in}\; G(\hat \Q^{\square})\}.$$
  This can be seen by the theory of log abelian varieties with constant degeneration.
   Then we have such an isomorphism (after taking localizations) at any point of some open neighborhood of $u$. 
\end{sbpara}

\section{Complements} 
\label{s:toroidal}

  Here we give some complements. 
 
  In Section \ref{ss:overC},
 we explain the relationship of our toroidal compactification with the analytic moduli space (\cite{KKN1}).

  In Sections \ref{s:intfan}--\ref{ss:Lan}, we explain the relationship of  our toroidal compactification with that of Lan \cite{Lan}.

Let $\cA_{\cH}$ be the moduli space of the moduli functor $F_{\cH}$ (\ref{FHbarH}) of abelian varieties with PEL structure. 
  It 
is denoted by ${\sf M}_{\cH}$ in \cite{Lan}. 

  Let $\bar \cA_{\cH, \Sig}$ be our moduli space of the moduli functor $\bar F_{\cH, \Sig}$ (\ref{FKSig2}) of log abelian varieties with PEL structure whose local monodromies are in $\Sig$ (cf.\ Theorem \ref{t:Main_Sig}).

\subsection{Complex analytic theory}
\label{ss:overC}

In \cite{Lan2}, Lan compares his algebraic toroidal compactification with the analytic toroidal compactification in \cite{AMRT}. 
By using our construction of the algebraic  toroidal compactification, we describe the relation with the analytic theory  by comparing log abelian varieties with analytic log abelian varieties in \cite{KKN1}. 

In this Section \ref{ss:overC}, we  take $\square =\emptyset$  as in \cite{Lan2}. 

The log complex analytic space $\bar \cA_{\cH, \Sig}^{\an}$ associated to $\bar \cA_{\cH, \Sig}\otimes_{F_0} \C$ is a toroidal compactification of the complex analytic space $\cA_\cH^{\an}$ associated to $\cA_\cH\otimes_{F_0} \C$, which is explicitly described as in Proposition \ref{overC} below. 
  We explain this by using the theory of analytic log abelian varieties in \cite{KKN1}.

\begin{sbpara}  We have
$$\cA_\cH^{\an}= \coprod_{\xi}   G^{(\xi)}(\Q) \bs (D \times G({\bf A}^f_\Q)/\cH)$$
(cf. \cite{GN} 3.6). The notation here is as follows.
Here $\xi$ ranges over all elements of the kernel of the map of Galois cohomologies $H^1(\Q, G)\to \prod_v H^1(\Q_v, G)$, where $v$ ranges over all places of $\Q$.  Each $\xi$ corresponds to the isomorphism class of a $B$-module $V^{(\xi)}$ with a skew-symmetric  form $\psi^{(\xi)} : V^{(\xi)}\times V^{(\xi)}\to \Q(1)$ such that  $(V^{(\xi)}, \psi^{(\xi)})_{\Q_v}\cong (V, \psi)_{\Q_v}$ 
for all $v$. We fix this isomorphism for each $v$. For each $\xi$, $G^{(\xi)}$ is the group of similitudes of $(V^{(\xi)}, \psi^{(\xi)})$.

$D$ is the conjugacy class of $h_0$ in $G(\R)$, which has a natural structure of a complex analytic manifold.   Since $G^{(\xi)}_\R= G_\R$, $G^{(\xi)}(\Q)$ acts on $D$ in the above presentation of $\cA_\cH^{\an}$. 

${\bf A}^f_\Q$ denotes the ring $\hat \Q^{\square}$  of finite adeles of $\Q$. 

For $(A, \iota, p, \eta)\in \cA_{\cH}(\C)=\cA_{\cH}^{\an}$, the corresponding element of the set on the right hand side is given as follows. Consider  $V_A:=H_1(A^{\an}, \Q)$  endowed with  a $B$-module structure given by $\iota$ and a skew-symmetric form  given by the polarization $p$.  By \cite{Ko} Lemma 4.2, we have an isomorphism $V_A\otimes \R\cong V \otimes \R$ of   $B$-modules endowed with such forms. By the level structure $\eta$, we have also  $V_A \otimes {\bf A}^f_\Q \cong V \otimes {\bf A}^f_\Q$. Hence the difference of $V$ and $V_A$ gives an element $\xi$ of $H^1(\Q, G)$ which becomes trivial in $H^1(\Q_v, \Q)$ for every place $v$ of $\Q$. Fix an isomorphism $V_A\cong V^{(\xi)}$. Then by \cite{Ko} Lemma 4.3, we have an element of $D$, and the level structure gives an element of $G({\bf A}^f_\Q)/\cH$. 
\end{sbpara}

\begin{sbpara} For each $\xi$, we define a rational fan $\Sig(\xi)$ in  $\Lie(G^{(\xi)})_\R$ as follows.

Let $\Sig(\xi)$  be the set of all cones $\sig$ which are obtained in the following way. Take a  
$B$-submodule $J$ of $V^{(\xi)}$ such that $J^{\perp} \subset J$, where $J^{\perp}$ is the annihilator of $J$ for $\psi^{(\xi)}$,  and let $W=V^{(\xi)}/J$. Let $g$ be the composite map $V\otimes {\bf A}^f_\Q\cong V^{(\xi)} \otimes {\bf A}^f_\Q \to W \otimes {\bf A}^f_\Q$. Let $\tau\in \Sig(W, g)$, and let $\sig= \{N_\mu\;|\; \mu\in \tau\}$ (\ref{cond4-0} (3)). 

Here for a symmetric bilinear form $\mu: W \times W \to \R$ such that $\mu(bx, y)=\mu(x, b^*y)$ for all $b\in B$ and $x, y\in W$, we have an element $N_{\mu}\in \Lie(G^{(\xi)})_\R$ by linearity because in the case $\mu(W \times W)\subset \Q$, 
$N_\mu$ in \ref{cond4-0} (3) satisfies $N_\mu(V^{(\xi)})\subset V^{(\xi)}$ and belongs to $\Lie(G^{(\xi)})$. 

We prove that $\Sig(\xi)$ is actually a fan. For this, it is sufficient to prove that for $\sig=\{N_{\mu}\;|\; \mu \in \tau\}$, $\sig'= \{N_{\mu}\;|\; \mu \in \tau'\}$ ($\tau\in \Sig(W, g)$, $\tau'\in \Sig(W', g')$), $\sig\cap \sig'$ is a face of $\sig$. 
 Take an interior point $N$ of the cone 
$\sig \cap \sig'$, 
write $N=N_{\mu}=N_{\mu'}$ ($\mu\in \tau$, $\mu'\in \tau'$), 
let $\tau_1$ (resp.\ $\tau'_1$) be the smallest face of $\tau$ (resp.\ $\tau'$) which contains $\mu$ (resp.\ $\mu'$), 
and let $J_1$ (resp.\ $J'_1$) be the inverse image of $\Ker(\mu)$ (resp.\ $\Ker(\mu')$) in $V^{(\xi)}$. 
Then $J_1=\{x\in V^{(\xi)}\;|\; N(x)=0\}=J'_1$. Let $W''=V^{(\xi)}/J_1= V^{(\xi)}/J_1'$. 
Since $\Ker(\mu)\subset \Ker(\nu)$ for all $\nu\in \tau_1$ and $\Ker(\mu') \subset \Ker(\nu)$ for all $\nu\in \tau'_1$,  $\tau_1$ (resp.\ $\tau'_1$) comes from some element $\tau_2$ (resp.\ $\tau'_2$) 
of $\Sig(W'',g'')$, where $g''$ is the map $V\otimes {\bf A}_{\Q}^f\cong V^{(\xi)}\otimes {\bf A}^f_{\Q} \to W'' \otimes {\bf A}^f_{\Q}$. Since  these $\tau_2$ and $\tau'_2$ share an interior point corresponding to $N$, we have $\tau_2=\tau'_2$. Hence $N$ belongs to $\alpha=\{N_{\nu}\;|\;\nu\in \tau_2=\tau'_2\}$ which is a face of $\sig$ and is a face of $\sig'$. 
  Since $\alpha$ is a face of $\sig$, it is a face of $\sig\cap \sig'$. 
  Since $\alpha$ contains an interior point $N$ of $\sig\cap \sig'$, we have $\alpha=\sig\cap \sig'$. 
  Hence  $\sig\cap \sig'$ is a face of $\sig$. 
\end{sbpara}

\begin{sbpara} For $g\in G({\bf A}^f_\Q)$, let $\Gamma(\xi, g)$ be the arithmetic subgroup $G^{(\xi)}(\Q) \cap g\cH g^{-1}$ of $G^{(\xi)}(\Q)$. 

By \cite{AMRT} and by
\cite{KKN1} Section 4, there is a toroidal compactification $\Gamma(\xi,g)\bs D_{\Sig(\xi)}$ of $\Gamma(\xi,g)\bs D$ associated to $\Sig(\xi)$ which is the moduli space of analytic log abelian varieties with PEL structure. 

\end{sbpara}

\begin{sbprop}\label{overC} Let $R\subset G({\bf A}^f_\Q)$ be a representative of $G^{(\xi)}(\Q)\bs G({\bf A}^f_\Q)/\cH$ (this $R$ is a finite set as is well-known). Then we have a commutative diagram

 $$\begin{matrix} \cA_\cH^{\an} & =& \coprod_{\xi} \coprod_{g\in R} \;\Gamma(\xi, g)\bs D\\
 \cap && \cap\\
\bar \cA_{\cH,\Sig}^{\an} &=&  \coprod_{\xi}\coprod_{g\in R}   \;\Gamma(\xi,g)\bs D_{\Sig(\xi)}.\end{matrix}$$

\end{sbprop}

\begin{pf} We prove this proposition by comparing log abelian varieties  with the analytic log abelian varieties in \cite{KKN1}. Let $Y= \coprod_{\xi}\coprod_{g\in R}   \;\Gamma(\xi,g)\bs D_{\Sig(\xi)}$. By  \cite{KKN2} Theorem 4.10,
for a log algebraic space $S$ in the first sense over $\C$  locally of finite type and for a log abelian variety $A$ over $S$, the inverse image $A^{\an}$ of $A$ on  the category of fs log analytic spaces over $S^{\an}$ is 
  an analytic log abelian variety over $S^{\an}$.  In particular, for the universal log abelian variety $A$ over $S=\bar \cA_{\cH, \Sig}$, $A^{\an}$ is an analytic log abelian variety over $S^{\an}$. 
Hence we have a canonical morphism from $S^{\an}= \bar \cA_{\cH, \Sig}^{\an}$ to the moduli space $Y$ of analytic log abelian varieties. 

It remains to prove that this morphism $\bar \cA_{\cH, \Sig}^{\an}\to Y$ is an isomorphism. For this, we use the fact that a morphism $Z_1\to Z_2$ of  fs log analytic spaces is an isomorphism if and only if the induced map  $\text{Mor}(S, Z_1)\to \text{Mor}(S, Z_2)$  of the sets of morphisms  is bijective  for every fs log analytic space $S$ whose underlying set is a one point set. Such an $S$ is written as

\medskip

$(*)$ $S=\Spec(R)$ for a local ring $R$ over $\C$ such that $R$ is of finite dimensional as a $\C$-vector space, endowed with an fs log structure

\medskip
\noindent
(we have $S^{\an}=S$ as a locally ringed space with a log structure). We prove that 
$\text{Mor}(S, \bar \cA^{\an}_{\cH, \Sig})\to \text{Mor}(S, Y)$ is bijective for $S$ as in $(*)$.  
  For such an $S$, we have $\text{Mor}(S, \bar \cA^{\an}_{\cH, \Sig})=\bar F_{\cH, \Sig}(S)$ and $\text{Mor}(S, Y)= \bar F_{\an, \cH, \Sig}(S)$, where $\bar F_{\an, \cH, \Sig}$ is the evident version of $\bar F_{\cH, \Sig}$ for fs log analytic spaces defined by replacing log abelian varieties by analytic log abelian varieties. Hence we are reduced to the fact that for such an $S$, the functor $A\mapsto A^{\an}$ from the category $\cC_1$ of polarized log abelian varieties over  $S$ to the category $\cC_2$  of polarized analytic log abelian varieties over $S$ is an equivalence of categories. We have equivalences of categories $$\cC_1 \simeq \cC_1',  \quad \cC_2\simeq \cC_2',$$
where $\cC_1'$ is the category of polarized log $1$-motifs (\cite{KKN2} Section 2) over $S$ and $\cC_2'$ is the category of polarized log Hodge structures on $S$ of weight $-1$ which belong to the category $\cH_S$ in  \cite{KKN1} 3.1.3 (these equivalences are by \cite{KKN2} Theorem 3.4 and by \cite{KKN1} Theorem 3.1.5, respectively). The converse functor  of $\cC_1\to \cC_2$ is induced by the functor $\cC_2'\to \cC_1'$ defined as follows. Let $H$ be an object of $\cC_2'$,
 Then the corresponding log $1$-motif is $[Y\to \cG_{\log}]\in \cC'_1$  given as follows.
  Let $\tau: S^{\log}\to S$ be the canonical proper surjective map. 
  Let $H_{\Z}$ be the $\Z$-structure of $H$ which is a locally constant sheaf on $S^{\log}$, and let $F$ be the Hodge filtration on the vector bundle $H_{\cO}$ on $S$.  We have an isomorphism $\cO_S^{\log} \otimes _\Z H_\Z\cong \cO^{\log}_S \otimes_{\tau^{-1}(\cO_S)} \tau^{-1}(H_{\cO})$ on $S^{\log}$. 
   We have the relative monodromy filtration  $W'$ on $H_{\Z}$ given as follows: $W'_wH_\Z=H_\Z$ for $w\geq 0$, $W'_{-1}H_\Z= \tau^{-1}\tau_*H_\Z$, $W'_{-2}H_\Z$ is 
    the annihilator of $W'_{-1}H_\Z$ for the polarization, and $W'_wH_\Z=0$ for $w\leq -3$. The sheaves $W'_{-1}H_\Z$, $H_\Z/W'_{-2}H_\Z$ and $\gr^{W'}_wH_\Z$ for $w\in \Z$ coincide with their $\tau^{-1}\tau_*$ and hence they are constant sheaves and regarded as sheaves on the point $S$. 
   Let  $Y=\gr^{W'}_0H_\Z$ regarded as a sheaf on $S$. Let $\cG$ be the semiabelian scheme over $S$ whose mixed Hodge structure $H_1(\cG)$ over $S$ is as follows. Its $\Z$-structure $H_1(\cG^{\an}, \Z)$ is $W'_{-1}H_\Z$, which we regard as a sheaf on $S$, and its  weight (resp.\ Hodge) filtration is  the restriction of $W'$ (resp.\ $F$). For a free $\Z$-module $\Lambda$ of finite rank which we regard as a Hodge structure of weight $0$ over $S$ in the natural way,  we have 
   $$\Ext(\Lambda, H_1(\cG))\cong \Hom(\Lambda, \cG), \quad  \Ext_{\log}(\Lambda, H_1(\cG))\cong \Hom(\Lambda, \cG_{\log}),$$
   where $\Ext$ (resp.\ $\Ext_{\log}$) is the group of extension classes in the category of mixed Hodge structures (resp.\ log mixed Hodge structures (\cite{KKN1} Section 2)) over $S$. These isomorphisms are obtained by the method of \cite{De2} (10.1.3) on $1$-motifs.  Take $\Lambda=Y$. Then the exact sequence $0\to H_1(\cG) \to (H_{\Z}, W', F) \to Y \to 0$ of log mixed Hodge structures over $S$ gives an element of  $\Ext_{\log}(Y, H_1(\cG))$ and hence gives a log $1$-motif $Y\to \cG_{\log}$. 
The polarization of this log $1$-motif is induced from the polarization of $H$. 
\end{pf}

\subsection{More integral formulation of fan}\label{s:intfan}

  We will explain that Lan's compactification is the underlying space of our toroidal compactification.
  
  For this comparison, first in this section, 
we rewrite the formulation of the fan in Section \ref{s:fan} in an 
equivalent form (cf.\ \ref{S3S4}), which is more integral and is nearer to the formulation of Lan.  

  In the following, we assume that $\cO$ and $L$ in Section \ref{s:integral} are given and we assume $\cH\subset G_L(\Z^{\square})$. 

\begin{sbpara}\label{intSig1}  A {\it more integral compatible family of complete fans} for $(B, V, \psi,\cH)$ is to give a complete fan $\Sig(W,g)$ for $W\otimes \Q$ (\ref{cofan})  for each pair $(W, g)$ of a finitely generated $\cO$-module  $W$ which is torsion-free as a $\Z$-module and a surjective $\cO\otimes \hat \Z^{\square}$-homomorphism $$g: L \otimes \hat \Z^{\square}\to W\otimes \hat \Z^{\square}$$ whose kernel $J$ satisfies
$J^{\perp} \subset J$ ($J^{\perp}\subset L\otimes\hat \Z^{\square}$ denotes the annihilator of $J$ under $\psi$), satisfying the following conditions (i) and (ii).

(i) $\Sig(W,gk)= \Sig(W,g)$ for $k\in \cH$.

(ii) Let $W'$ be a finitely generated $\cO$-module which is torsion-free as a $\Z$-module and  let $\gamma: W\to W'$ be a surjective $\cO$-homomorphism. Then $\Sig(W',\gamma g)= \gamma \Sig(W,g)$. Here $\gamma g$ denotes the composite homomorphism $$L\otimes \hat \Z^{\square} \overset{g}\to W\otimes \hat \Z^{\square}\overset{\gamma}\to W'\otimes \hat \Z^{\square},$$ 
and $\gamma \Sig(W,g)$ is the set of cones $\sig$ in $C^+(W'\otimes \Q)$ such that the image of $\sig$ in $C^+(W\otimes \Q)$ under the map $C^+(W'\otimes \Q) \to C^+(W\otimes \Q)$ induced by $\gamma$ belongs to $\Sig(W,g)$. 
\end{sbpara}

\begin{sbpara}\label{S3S4}  $\Sig$ in Section \ref{s:fan} gives a more integral $\Sig$ in this section. 
  For a $W$ in \ref{intSig1}, 
$W \otimes \Z_{(\square)}$ is a $W$ in Section \ref{s:fan}. 
  Hence $\Sig$ in Section \ref{s:fan} gives a $\Sig$ in \ref{intSig1}.

  Conversely, for a $W$ in Section \ref{s:fan}, the intersection of the image of $L \otimes \hat \Z^{\square} \to W\otimes \hat \Q^{\square} $ and $W$ is a $W$ in \ref{intSig1}. 
  Hence $\Sig$ in \ref{intSig1} gives a $\Sig$ in Section \ref{s:fan}.  

  In the following, we identify both $\Sig$. 
\end{sbpara}

\begin{sbpara}\label{intSig2}
 Let $\Sig=(\Sig(W,g))_{(W,g)}$ be a more integral compatible family of complete fans for $(B, V,\psi, \cH)$.

We can rewrite the moduli functor $$\bar F_\cH\supset \bar F_{\cH, \Sig}: (\fs/O_{F_0, (\square)})\to (\mathrm{set})$$
as that $\bar F_{\cH, \Sig}(S)$ is the part of $\bar F_\cH(S)$ consisting of elements satisfying the following condition (i).
 
 \medskip

(i)  For each geometric point $s$ of $S$, let $W={\bar Y}_s$, and let $g: L \otimes \hat \Z^{\square}\to W \otimes \hat \Z^{\square}$ be the surjective homomorphism induced by the inverse of the level structure of $A$. 
  Then 
there exists a $\sig\in \Sig(W, g)$ such that for every $a\in \Hom((M_S/\cO^\times_S)_s,\N)$, the composition $W \times W \to  (M_S^{\gp}/\cO_S^\times)_s \overset{a}\to \bR$ belongs to $\sig$.
\end{sbpara}

\begin{sbpara}\label{intSig4}  
We can also rewrite the moduli functor for a single $\sigma$ as follows. 

In \ref{intSig1}, fix a $(W,g)$ and a $\sig \in \Sig(W,g)$. 
  Then $\bar F'_{\cH, \sig}(S)$ is identified with the set of pairs of (the same thing as in the definition of $\bar F_\cH$) and (a surjective $\cO$-homomorphism $g': W\to \bar Y$) such that the composition $L\otimes \hat \Z^{\square}\overset{g}\to W \otimes \hat \Z^{\square}\overset{g'}\to \bar Y\otimes \hat \Z^{\square}$ is induced by the inverse of the level structure and such that for every geometric point $s$ of $S$ and every $a\in \Hom((M_S/\cO^\times_S)_s,\N)$, the bilinear form 
$W \times W \to \bR$ induced by $g'$ and $a$ belongs to $\sig$.
\end{sbpara}

\subsection{Comparison with Lan's compactification}
\label{ss:Lan}
In this section we explain
that Lan's compactification is the underlying space of our space. 
(As explained in Remark \ref{smcone}, we allow slightly  more general  fans  than the ones in \cite{Lan}.)
   In the following, we use the notation in Lan's \cite{Lan} freely. 
  The assumptions are the same as in \cite{Lan}.  In particular, we assume here \cite{Lan} Condition 1.4.3.10: 

\smallskip

\noindent $(*)$ The action of $\cO$ on $L$ extends to an action of some maximal order in $B$ containing $\cO$. 

\smallskip

\begin{sbpara}
  Take a compatible choice $\Sig'=\{\Sig'_{\Phi_{\cH}}\}_{[(\Phi_{\cH},\delta_{\cH})]}$ of admissible smooth rational polyhedral cone decomposition data in the sense of \cite{Lan} Definition 6.3.3.4 for ${\textsf M}_{\cH}$. 
  We compare Lan's compactification ${\textsf M}_{\cH,\Sigma'}^{\mathrm{tor}}$ with our space $\bar \cA_{\cH,\Sig}$, where $\Sig$ is determined by $\Sig'$ as follows. 
\end{sbpara}

\begin{sbpara}
  We explain how to associate to $\Sig'$ a more integral compatible family of complete fans $\Sig$ in Section \ref{s:intfan}.
  Recall that $\Sig'$ is indexed by the cusp labels at level $\cH$ in the sense of \cite{Lan} Definition 5.4.2.4, that our $\Sig$ is indexed by the set of $(W,g)$, and that both index sets relate to the stratifications of the respective compactifications. 
  We would like to associate to each $(W,g)$ as in \ref{intSig1} a cusp label at level $\cH$.  
  Recall that a cusp label at level $\cH$ is an equivalence class of triples $({\tt Z}_{\cH}, \Phi_{\cH}, \delta_{\cH})$ of some linear algebraic data. 
  
  Let $(W,g)$ be as in \ref{intSig1}.
  We assume that the following condition is satisfied.

\smallskip

\noindent 
$(**)$ 
$g(L)$ is a subgroup of $W$ of finite index. 

\smallskip

  We associate to $(W,g)$ a cusp label at level $\cH$ under this condition. 

  First we give 
a collection ${\tt Z}_{\cH}=\{{\tt Z}_{H_n}\}_n$ 
of $H_n$-orbits of fully symplectic-liftable admissible filtrations on $L/nL$ with respect to $(L,\langle \cdot,\cdot\rangle)$. 
  Let $n$ be a positive integer such that $\square\nmid n$ and ${\cal U}^{\square}(n) \subset \cH$. 
  Consider the 3-step filtration ${\tt Z}_n$ on $L/nL$ determined by the images of $J$ and $J^{\perp}$ in \ref{intSig1}. 
  Then ${\tt Z}_n$ determines an $H_n$-orbit ${\tt Z}_{H_n}$ of admissible filtrations on $L/nL$ that are fully symplectic-liftable. 

  Next we give 
a torus argument $\Phi_{\cH}=\{\Phi_{H_n}\}_n$ at level $\cH$ for ${\tt Z}_{\cH}$. 
  Let $X=\Hom_{\bZ}(J^{\perp} \cap L,\bZ)$ and $Y=g(L)$. 
  Then $g$ and $\psi$ induce 
an $\cO$-linear injection 
$\phi: Y \to X$. 
  For each $n$ as above, 
  let $\varphi_{-2,n}:
\mathrm{Gr}_{-2,n}^{{\tt Z}_n} \overset \sim \to \Hom_{\hat \bZ{}^{\square}}(X/nX, \bZ/n\Z(1))$ and 
  $\varphi_{0,n}:
\mathrm{Gr}_{0,n}^{{\tt Z}_n} \overset \sim \to Y/nY$ 
be the natural isomorphisms. 
  Let $\Phi_n=(X, Y, \phi, \varphi_{-2,n}, \varphi_{0,n})$, which is a torus argument at level $n$ for ${\tt Z}_n$. 
  Let $\Phi_{H_n}$ be the $H_n$-orbit of $\Phi_n$. 
  Then $\Phi_{\cH}=\{\Phi_{H_n}\}_n$ is a torus argument at level $\cH$ for ${\tt Z}_{\cH}$. 

  Third, by the condition $(*)$ in the above, 
there is a liftable splitting $\delta_{\cH}$ of ${\tt Z}_{\cH}$, and the triple $({\tt Z}_{\cH}, \Phi_{\cH}, \delta_{\cH})$ determines a cusp label at level $\cH$, which is the cusp label corresponding to $(W,g)$. 

  Let $\Sig(W,g)$ be the complete fan induced by $\Sig'_{\Phi_{\cH}}$. 

  For a general $(W,g)$, take a $k \in \cH$ such that $(W,gk)$ satisfies the condition $(**)$ and 
let $\Sig(W,g)$ be the one induced by $\Sig(W,gk)$.
  Let $\Sig=(\Sig(W,g))_{(W,g)}$. 
\end{sbpara}

\begin{sbrem}\label{smcone}
  Note that we treat essentially more $\Sig=(\Sig(W,g))_{(W,g)}$ than the $\Sig'$s treated in \cite{Lan}: 
  For technical reasons, \cite{Lan} imposes on $\Sig'$ the condition of the smoothness and the condition for avoiding self-intersections in the construction. 
  Thus a more integral compatible family $\Sig$ of complete fans can be obtained from some $\Sig'$ in the above way only if all cones in all $\Sig(W,g)$ are smooth and only if each $\Sig(W,g)$ satisfies the condition corresponding to Condition 6.2.5.25 in \cite{Lan}.
  Recall that in the framework of log geometry any fan behaves as if it were a smooth fan, which is the reason why the smoothness of fans is not essential to us.  
\end{sbrem}

\begin{sbprop}
\label{p:comparison}
  Let $\Sig'$ and $\Sig$ be as above. 
  Then the space $\bar \cA_{\cH,\Sig}$ representing the functor 
$\bar F_{\cH,\Sig}$ in Theorem $\ref{t:Main_Sig}$ 
coincides with Lan's compactification ${\mathrm{\sf M}}_{\cH,\Sigma'}^{\mathrm{tor}}$ 
({\rm{\cite{Lan}}} Theorem $6.4.1.1$) 
of the moduli space ${\sf M}_{\cH}$ of abelian varieties with PEL structure, 
endowed with the fs log structure defined by the divisor ${\sf M}_{\cH,\Sigma'}^{\mathrm{tor}}-{\sf M}_{\cH}$. 
\end{sbprop}

\begin{sbpara}
  We sketch the proof, which is parallel to the non-coefficient case in \cite{KKN7} Section 9.

  We construct a universal family of log abelian varieties with coefficients over Lan's compactification. 
  Let $L={\textsf M}_{\cH,\Sig'}^{\mathrm{tor}}$ be Lan's compactification endowed with the log structure defined by ${\textsf M}_{\cH,\Sig'}^{\mathrm{tor}}-{\textsf M}_{\cH}$. 
  Consider the natural morphism from $L$ to the space $L_0$ without coefficients, that is, the space of Faltings--Chai endowed with the log structure.
  By \cite{KKN7} Section 9, $L_0$ coincides with the fine moduli of log abelian varieties without coefficients, and 
the family of Faltings--Chai is regarded as a model of  
the universal family of log abelian varieties over $L_0$.
  By pulling back this universal family, we have a family $A$ of log abelian varieties on $L$ one of whose models is the family of degenerating abelian varieties of Lan. 
  The last family has a coefficient $\cO \to \End(\cG)$, where $\cG$ is the semiabelian part of $A$. 
  We see that it factors through $\End(A) \subset \End(\cG)$ by checking formally so that it gives a coefficient $\cO \to \End(A)$ on $A$.  
  Thus we have a family of log abelian varieties with coefficients over $L$.

  Since our space is a fine moduli, it gives a morphism from $L$ to our space 
$\bar \cA_{\cH,\Sig}$. 
  We can see that this morphism is formally an isomorphism (so that it is an isomorphism) by considering the local moduli as in \cite{KKN6} (cf.\ \ref{local_moduli}) and comparing it with the description of the completions of Lan's compactification in \cite{Lan} 6.2. 

\end{sbpara}

\noindent {\bf Correction to \cite{KKN5}.}
  In the end of \cite{KKN7}, we corrected a part of \cite{KKN5}, that is, 
we suppressed \cite{KKN5} Proposition 12.8 (4), \cite{KKN5} Remark 12.8.1 (1), \cite{KKN5} Lemma 12.9, and \cite{KKN5} Proposition 12.11. 
  Here we give a recovery for the cases of $\Gmlog$ of them. 
  Let $A$ be a weak log abelian variety over a locally noetherian fs log scheme $S$. 
  Then, since $A$ is covered by representable objects, the log Hilbert 90 (\cite{Kato:FI2} Corollary 5.2) implies that 
$H^1_{\et}(A,\Gmlog)=H^1_{\ket}(A,\Gmlog)$.
  By this, the cubic isomorphism for $\Gmlog$-torsors on the k\'et site under the assumption of \cite{KKN5} Theorem 2.2 (c) is reduced to that on the \'etale site. 
  Hence we recover the missing exact sequence in the proof of \cite{KKN5} Proposition 12.8 (4). 
  Therefore we recover the conclusions of the case of $F=\Gmlog$ of Proposition 12.8 (4), the case of $\Gmlog$ of \cite{KKN5} Remark 12.8.1 (1), and the case $F'=\Gmlog$ of \cite{KKN5} Lemma 12.9.

\noindent Takeshi Kajiwara

\noindent Department of Applied mathematics \\
Faculty of Engineering \\
Yokohama National University \\
Hodogaya-ku, Yokohama 240-8501 \\
Japan

\noindent kajiwara-takeshi-rj@ynu.ac.jp 
\par\bigskip\par

\noindent Kazuya Kato

\noindent 
Department of Mathematics
\\
University of Chicago
\\
5734 S.\ University Avenue
\\
Chicago, Illinois, 60637 \\
USA

\noindent kkato@math.uchicago.edu
\par\bigskip\par

\noindent Chikara Nakayama

\noindent Department of Economics \\ Hitotsubashi University \\
2-1 Naka, Kunitachi, Tokyo 186-8601 \\ Japan

\noindent c.nakayama@r.hit-u.ac.jp
\end{document}